\documentclass[a4paper]{article}
\usepackage[utf8x]{inputenc}
\usepackage{enumerate}
\usepackage{amssymb,amsmath,amsthm,amscd,latexsym}
\usepackage{setspace}
\usepackage{graphics, graphicx}
\usepackage{color}

\newtheorem{theo}{Theorem}[section]

\newtheorem{corol}[theo]{Corollary}
\newtheorem{prop}[theo]{Proposition}

\newtheorem{lem}[theo]{Lemma}
\newtheorem{slem}[theo]{Sublemma}
\newcounter{tmp}

\theoremstyle{definition}
\newtheorem{defi}[theo]{Definition}
\newtheorem{rem}[theo]{Remark}
\newtheorem*{example}{Example}

\newcommand{\can}{\mathrm{can}}
\newcommand{\diff}{\mathrm{Diff}}
\newcommand{\inj}{\mathrm{inj}}

\newcommand{\Rm}{\mathrm{Rm}}
\newcommand{\Ric}{\mathrm{Ric}}

\newcommand{\Rmin}{\mathrm{R_{min}}}

\newcommand{\NN}{\mathbf{N}}
\let\nn\NN
\newcommand{\Int}{\mathrm{Int}}

\newcommand{\std}{\mathrm{std}}
\newcommand{\RR}{\mathbf{R}}
\let\rr\RR

\newcommand{\cR}{\mathcal{R}}
\newcommand{\cF}{\mathcal{F}}
\newcommand{\cS}{\mathcal{S}}
\newcommand{\cA}{\mathcal{A}}
\newcommand{\cB}{\mathcal{B}}
\newcommand{\cC}{\mathcal{C}}
\newcommand{\cX}{\mathcal{X}}
\newcommand{\cP}{\mathcal{P}}
\newcommand{\cQ}{\mathcal{Q}}

\newcommand{\gcyl}{g_{\mathrm{cyl}}}
\newcommand{\surg}{\mathrm{surg}}

\let\bord\partial
\let\bydef\emph
\let\epsi\varepsilon

\newcommand{\Addresses}{{
  \bigskip
  \footnotesize

 \textsc{Institut de Mathématiques de Bordeaux, Université de Bordeaux,
    33405, Talence, France}\par\nopagebreak
  \textit{E-mail address}: \texttt{laurent.bessieres@math.u-bordeaux.fr}

  \medskip

 \textsc{Institut Fourier, Université Grenoble Alpes, 100 rue des maths, 38610 Gières, France}\par\nopagebreak
  \textit{E-mail address}: \texttt{g.besson@univ-grenoble-alpes.fr}

  \medskip

 \textsc{Institut Montpelli\'erain Alexander Grothendieck, CNRS, Université de Montpellier}\par\nopagebreak
  \textit{E-mail address}: \texttt{sylvain.maillot@univ-montp2.fr}
 \medskip

 \textsc{Princeton University, Fine Hall, Princeton NJ 08544, USA}\par\nopagebreak
  \textit{E-mail address}: \texttt{coda@math.princeton.edu}

}}

\title{Deforming $3$-manifolds of bounded geometry and uniformly positive scalar curvature}
\author{Laurent Bessi\`eres, G\'erard Besson, Sylvain Maillot\\ and Fernando C. Marques}

\begin{document}
\begin{spacing}{1.2}

\maketitle

\begin{abstract}
We prove that the moduli space of complete Riemannian metrics of bounded geometry and uniformly positive scalar curvature on an orientable 3-manifold is path-connected. This generalizes the main result of the fourth author~\cite{marques:deforming} in the compact case. The proof uses Ricci flow with surgery as well as arguments involving performing infinite connected sums with control on the geometry.
\end{abstract}

\section{Introduction}\label{sec:intro}

\subsection{Summary of earlier results}
A classical topic in Riemannian geometry is the study of complete manifolds
with positive scalar curvature. A first basic problem is that of existence: in a given dimension, characterize the class of smooth manifolds which admit such a metric. A second one is, given a manifold $M$ in this class, to determine whether such a metric is unique up to isotopy---possibly working modulo the action of the diffeomorphism group of $M$ on the space of metrics under consideration.

When the manifold is non-compact, each of these two problems gives rise to several questions, according to whether one studies metrics with positive scalar curvature or \bydef{uniformly positive} scalar curvature, i.e.~scalar curvature greater than some positive constant, or restricts attention to metrics satisfying additional bounds on the geometry.

In this paper we are concerned with dimension 3. We first briefly review known results in other dimensions. The Gauss-Bonnet formula implies that $S^2$ is the only two-dimensional orientable compact manifold with positive scalar curvature. Weyl~\cite{Weyl} proved that the space of such metrics on $S^2$ is path-connected, as a consequence of the Uniformization Theorem (later Rosenberg-Stolz~\cite{RosenbergStolz} proved contractibility).  In high dimensions this is no longer true. The space of metrics of positive scalar curvature on $S^n$ is disconnected for all dimensions $n\geq 7$ of the form $n=8k, n=8k+1$ or $n=4k-1$ (Hitchin~\cite{Hitchin}, Carr~\cite{carr}, Kreck-Stolz~\cite{KreckStolz}). We would also like to mention that similar questions have been studied in dimension~2 in a noncompact setting, see~Belegradek-Hu~\cite{bh:connectedness,bh:erratum}.

Let $M$ be a closed, orientable, connected, smooth 3-dimensional manifold.  The question of when such a manifold admits a metric of positive scalar curvature has been studied by Gromov-Lawson~\cite{Gro-Law:scalar, GromovLawson} and
Schoen-Yau~\cite{SchoenYau} in the late 1970s and completely settled by Perelman in his celebrated work on the Geometrisation Conjecture in 2003. Obvious examples of manifolds in this class are $S^2\times S^1$ and \emph{spherical manifolds} (i.e.~manifolds that are diffeomorphic to metric quotients of the round $3$-sphere.) Gromov-Lawson and Schoen-Yau observed independently that this class is closed under connected sum. Perelman completed the classification:
\begin{theo}[Perelman~\cite{Per2}]\label{thm:classif compact}
A closed, orientable, connected 3-manifold admits a Riemannian metric of positive scalar curvature if and only if it is a connected sum of spherical manifolds and $S^2\times S^1$'s.
\end{theo} 

We now state the uniqueness result due to the fourth author. For this we introduce some notation. Let $\sigma$ be a positive constant. Let $\mathcal{R}_+(M)$ (resp.~$\mathcal{R}_\sigma(M)$) denote the set of all metrics on $M$ which have positive scalar curvature (resp. scalar curvature $\ge \sigma$.) 

\begin{theo}[Marques~\cite{marques:deforming}]\label{thm:marques deforming}
Let $M$ be a closed, orientable, connected 3-manifold such that $\mathcal{R}_+(M)\not=\emptyset$. Then $\mathcal{R}_+(M)/\diff(M)$ is path-connected in the $\mathcal{C}^\infty$ topology.   
\end{theo}

The main theorem of this article is a generalisation of Theorem~\ref{thm:marques deforming} to possibly non-compact manifolds. Before we state it, we review some results on the existence question in the setting of open 3-manifolds.

Let $M$ be an open, orientable, connected 3-manifold. The basic question of when does such a manifold admit a complete metric $g$ of positive scalar curvature (\cite{Yau}, Problem 27) is still wide open. Progress has been made under the stronger hypothesis of  uniformly positive scalar curvature. In particular, Cheng, Weinberger and 
Yu~\cite{cwy:taming} classified complete 3-manifolds with uniformly positive scalar curvature and finitely generated fundamental group. More closely related to this paper is a result of the first three authors~\cite{B2M:scalar} which we now state.

Recall that a Riemannian manifold is said to have \emph{bounded geometry} if it has bounded sectional curvature and positive injectivity radius.  We also need to recall the definition of a connected sum along a locally finite graph. For the moment, we will assume that the summands are closed---later we will need to consider non-compact summands (cf.~Section~\ref{sec:GL}.)

Let $\cX$ be a collection of closed, oriented $3$-manifolds and $G$ be a locally finite connected graph. Fix a map $X:v \mapsto X_v$ which associates to each vertex of $G$ a copy of some manifold in $\cX$. To the pair $(G,X)$ we associate a $3$-manifold $M$ in the following way: for each $v$ we let $Y_v$ be $X_v$ with $d_v$ punctures, where $d_v$ is the degree of $v$. Let $Y$ be the disjoint union of all $Y_v$'s. Then $M$ is obtained from $Y$ by the following operation: for each edge $e$ of $G$, call $v,v'$ the vertices connected by $e$; choose a $2$-sphere $S\subset \bord Y_v$ and a $2$-sphere $S'\subset \bord Y_{v'}$, and glue $Y_v$ and $Y_{v'}$ together along $S$ and $S'$ using an orientation-reversing diffeomorphism. We say that a $3$-manifold is a \emph{connected sum of members of $\cX$} if it is diffeomorphic to a manifold obtained by this construction.

When $G$ is a finite tree, we recover the usual notion of finite connected sum. When $G$ is a finite graph, we can turn it into a tree at the expense of adding some $S^2\times S^1$ factors. In the reverse direction, one can add edges to the graph and remove $S^2\times S^1$ factors. Thus the conclusion of Theorem~\ref{thm:classif compact} can be reformulated as: $M$ is a connected sum of spherical manifolds along a finite graph.

The generalisation of Theorem~\ref{thm:classif compact} due to the first three authors reads:
\begin{theo}[Bessi\`eres-Besson-Maillot~\cite{B2M:scalar}]\label{thm:BBM11}
Let $M$ be an orientable, connected $3$-manifold. Then $M$ admits a complete Riemannian metric of uniformly positive scalar curvature and bounded geometry if and only if there exists a finite collection $\cF$ of spherical $3$-manifolds such that $M$ is a connected sum of members of $\cF$.
\end{theo}

An explanation on how to deduce Theorem~\ref{thm:BBM11} from the main result of~\cite{B2M:scalar} is given in Appendix~\ref{sec:appendix}.

\subsection{Main results of this article}

{\em Throughout this paper we make the following convention: all manifolds are assumed to be smooth, orientable and without boundary. We will not assume that they are connected in general.}

\bigskip

Let $M$ be a $3$-manifold.  Given $\sigma>0$, we denote by $\mathcal{R}_\sigma^{bg}(M)$ the space of complete Riemannian metrics 
$g$ on $M$ with scalar curvature $\geq \sigma$ and bounded geometry.  We further set $\mathcal{R}_+^{bg}(M) = \bigcup_{\sigma>0} \mathcal{R}_\sigma^{bg}(M)$.
We endow $\mathcal{R}_+^{bg}(M)$ with the $C^\infty_{loc}$-topology. When $M$ is compact, every metric has bounded geometry, and we have $\mathcal{R}_+(M)=\mathcal{R}_+^{bg}(M)$ and $\mathcal{R}_\sigma(M)=\mathcal{R}_\sigma^{bg}(M)$ for each $\sigma$. 

The following is the main theorem of this article:
\begin{theo}\label{thm:connectedness}
Let $M$ be a connected $3$-manifold such that $\mathcal{R}_1^{bg}(M)\not=\emptyset$. Then $\mathcal{R}_1^{bg}(M)/\diff(M)$ is path-connected.   
\end{theo}

Since any metric in $\mathcal{R}_+^{bg}(M)$ differs from a metric in $\mathcal{R}_1^{bg}(M)$ by scaling, we obtain the following corollary:

\begin{corol}\label{corol:connectedness}
Let $M$ be a connected $3$-manifold such that $\mathcal{R}_+^{bg}(M)\not=\emptyset$. Then $\mathcal{R}_+^{bg}(M)/\diff(M)$ is path-connected.   
\end{corol}
 
\subsection{Strategy of proof}
Here we explain the main difficulty in extending Theorem~\ref{thm:marques deforming} to Theorem~\ref{thm:connectedness}. First we sketch the argument used in  \cite{marques:deforming} for proving Theorem~\ref{thm:marques deforming}.

Let $M$ be a closed, connected $3$-manifold such that $\cR_1(M)$ is nonempty.
By Theorem~\ref{thm:classif compact}, $M$ is a connected sum of spherical manifolds and $S^2\times S^1$'s. A crucial notion is that of a \emph{canonical metric} on $M$. To construct it, start with the unit sphere $S^3 \subset \RR^4$. For each spherical summand $\Sigma_i$, fix a round metric on $\Sigma_i$, a pair of points $(p_i^+,p_i^-)\in S^3\times \Sigma_i$ and perform a \emph{GL-sum} at the points $\{p_i^\pm\}$. This notion will be explained in  detail (as well as greater generality) in Section~\ref{sec:GL} below. Here we describe it informally: one removes a small metric ball $B_i^\pm$ around each $p_i^\pm$ and glues a thin tube between $\bord B_i^+$ and $\bord B_i^-$, preserving the condition of positive scalar curvature. Finally, for each $S^2\times S^1$ summand, do a `self GL-sum' of $S^3$, i.e.~a similar operation with both points $p_i^\pm$ in $S^3$.

Making this construction more precise and invoking Milnor's uniqueness theorem for the prime decomposition, one sees that the space of canonical metrics on $M$ is path-connected modulo diffeomorphism. Thus, for the purpose of proving Theorem~\ref{thm:marques deforming}, one can speak of `the' canonical metric on $M$.

The bulk of the proof consists in showing that any metric in $\cR_+(M)$ can be isotoped to the canonical metric. This uses a refined version of Perelman's Ricci flow with surgery (based on the monograph by Morgan and Tian~\cite{Mor-Tia}) as well as conformal deformations.

In the non-compact case, we can use Theorem~\ref{thm:BBM11} to recognize the topology of the manifold $M$. However, there is no uniqueness theorem for a presentation of $M$ as a connected sum of spherical manifolds. Thus as soon as the topology of $M$ becomes intricate, it is unclear how to define a canonical metric.

We shall define a special subset of $\cR_1^{bg}(M)$, consisting of what we call \emph{GL-metrics} (see below.) The proof has two parts. In the first part, we show that any two GL-metrics can be connected to each other (modulo diffeomorphism.) In the second part, we show that any metric in $\cR_1^{bg}(M)$ can be deformed into a GL-metric. We use arguments similar to those of~\cite{marques:deforming}, based on the version of Ricci flow with surgery developed in~\cite{B2M:scalar}. 

\subsection{Main technical results}
To make the above discussion more precise, we need to define GL-metrics. In order to do this, we need a topological notion, that of a \bydef{spherical splitting} of a 3-manifold, and the geometric notion of \bydef{straightness} with respect to a Riemannian metric.

\begin{defi}
Let $M$ be a $3$-manifold. A \emph{spherical system} in $M$ is a (possibly empty) locally finite collection $\cS$ of pairwise disjoint embedded $2$-spheres in $M$. The members of $\cS$ are called its \emph{components}. We denote by $M\setminus\cS$ the complement in $M$ of the union of all components of $\cS$.

A spherical system $\cS$ in $M$ is called a \emph{spherical splitting} if each connected component of $M\setminus\cS$ is relatively compact in $M$.
\end{defi}

In particular, the empty set is a spherical splitting if and only if $M$ is compact.

Let $g$ be a Riemannian metric on $M$. An embedded $2$-sphere $S\subset M$ is called \bydef{straight with respect to $g$} (or just \bydef{straight} if the metric is understood), if $S$ has an open tubular neighbourhood $U$, called a \emph{straight tube}, such that there is an isometry $f:U_0\to U$, where $U_0$ is a Riemannian product of a round 2-sphere with an interval, and the preimage of $S$  has the form $S^2 \times \{*\}$.  A spherical system is called \bydef{straight with respect to $g$} if all of its components are.

\begin{defi}\label{defi:GL metric}
Let $M$ be a $3$-manifold. A Riemannian metric $g\in\cR_+^{bg}(M)$ is called a \emph{GL-metric} if $M$ admits a spherical splitting which is straight with respect to $g$.
\end{defi}

We say that two metrics $g,g'\in\cR_1^{bg}(M)$ are \emph{isotopic} if there exists a continuous path $g_t \in \mathcal{R}_{1}^{bg}(M)$, for  $t\in [0,1]$, such that $g_0=g$ and $g_1= g'$. We say that $g,g'$ are \emph{isotopic modulo diffeomorphism} if $g$ is isotopic to a metric that is isometric to $g'$. 
One of the main technical results of this paper is the following:

\begingroup
\setcounter{tmp}{\value{theo}}
\setcounter{theo}{0} 
\renewcommand\thetheo{\Alph{theo}}
\begin{theo}\label{thm:gl connected}
Let $M$ be a 3-manifold and $g,g'\in\cR_1^{bg}(M)$ be two GL-metrics. Then $g$ is isotopic to $g'$ modulo diffeomorphism. 
\end{theo}

In order to state our second main technical result, we introduce more notation and terminology.  For $C>0$ we say that a metric on $M$ has \emph{geometry bounded by} $C$ if the norm of the curvature tensor is $\leq C$ and the injectivity radius is $ \geq C^{-1/2}$. Given $\sigma,C>0$ we denote by $ \mathcal{R}_\sigma^{C}(M)$ the set of complete Riemannian metrics on $M$ which have scalar curvature $\ge\sigma$ and geometry bounded by $C$. Two metrics $g,g' \in \mathcal{R}_{\sigma}^{C}(M)$ are said to be  \emph{isotopic in $\cR_1^C(M)$} if there exists a continuous path $g_t \in \mathcal{R}_{1}^{C}(M)$, $t\in [0,1]$, such that $g_0=g$ and $g_1= g'$. 

\begin{theo}\label{thm:isotopy_to_GL}
For every $A>0$ there exists $B=B(A)>0$ such that for every $3$-manifold $M$ and every metric $g \in \mathcal{R}_1^A(M)$, there exists a GL-metric $g' \in \mathcal{R}_1^{B}(M)$ which is isotopic to $g$ in $\mathcal{R}_1^{B}(M)$.
\end{theo}
\endgroup
\setcounter{theo}{\thetmp} 

It is immediate that Theorems~\ref{thm:gl connected} and~\ref{thm:isotopy_to_GL} together imply Theorem~\ref{thm:connectedness}.

\begin{rem}
Theorem~\ref{thm:isotopy_to_GL} is stronger than what is needed to prove Theorem~\ref{thm:connectedness} since we have a uniform bound on the geometry. If Theorem~\ref{thm:gl connected} could be similarly improved, i.e.~if one could show that for each $A>0$ there is $B>0$ such that for any 3-manifold $M$, any two GL-metrics in $\cR_1^A(M)$ are isotopic in $\cR_1^B(M)$, then one would get the stronger conclusion that for every $A>0$ there is a $B>0$ such that the map $\cR_1^A(M)/\diff(M) \to \cR_1^B(M)/\diff(M)$ induced by inclusion is 0-connected.
\end{rem}

Let us denote by $\diff^+(M)$ the group of orientation-preserving diffeomorphisms of $M$. One of the main ingredients for proving Theorem~\ref{thm:gl connected} is a version of the main result of~\cite{marques:deforming} which is slightly different from Theorem~\ref{thm:marques deforming}. We state it below for future reference.
\begin{theo}\label{thm:marques deforming2}
Let $M$ be a closed, oriented 3-manifold such that $\mathcal{R}_1(M)\not=\emptyset$. Then $\mathcal{R}_1(M)/\diff^+(M)$ is path-connected  in the $\mathcal{C}^\infty$ topology.  
\end{theo}

An explanation of how to deduce this from the arguments in~\cite{marques:deforming} is given in Appendix~\ref{appendix:marques}.

\medskip

The proof of Theorem \ref{thm:gl connected} relies on cut-and-paste arguments, which in our context means surgeries and GL-sums. Let us sketch this proof. Fix an orientation of $M$ and consider two GL-metrics $g,g'$ on $M$. Suppose first that $g,g'$ have a common straight splitting $\cS=\{S_\alpha\}$, and that $g$ coincides with $g'$ on a neighbourhood of $\cS$. Perform metric surgeries along $\cS$, cutting and gluing back standard caps. The resulting manifold $M_\#$ carries two metrics $g_\#$ and $g'_\#$ of positive scalar curvature coming from $g$ and $g'$ respectively. All the connected components of $M_\#$ are closed.  Theorem~\ref{thm:marques deforming2} applied to each component of $M_\#$ gives an isotopy modulo positive diffeomorphisms between $g_\#$ and $g'_\#$. Although we do not \emph{a priori} have a uniform upper bound on the geometry along this path of metrics, a reparametrisation trick (cf.~Subsection~\ref{subsec:same cS isotopic}) allows to perform GL-sums so as to obtain a path in $\cR_{1}^{bg}(M)$, connecting $g$ to $g'$ modulo diffeomorphism. In the general case of two arbitrary GL-metrics $g,g'$, the idea is to deform one of them so as to obtain two metrics which coincide near a straight splitting, as above. The proof is also by reduction to the compact case.  

The proof of Theorem~\ref{thm:isotopy_to_GL} uses \emph{surgical solutions}, a version of Ricci flow with surgery developed in \cite{B2M:scalar}. Starting from a Riemannian manifold $(M_0,g_0)$ with bounded geometry, it produces a sequence $(M_k,g_k(t))_{t \in [t_k,t_{k+1}]}$ of Ricci flows with uniformly bounded geometry, 
the jump from one interval to the next being realised by performing metric surgeries and discarding components with `nice' metrics. When the scalar curvature of the initial metric satifies $R_{g_0} \geq 1$,  the solution becomes extinct in finite time ($M_k=\emptyset$ for some $k$), and the last metric before extinction is nice enough so that it is straightforward to isotope it to some GL-metric. The construction of the isotopy to $g_0$ then uses a backward induction argument similar to the one
of~\cite{marques:deforming}, performing GL-sums of paths to GL-metrics. A technical issue here is to realise this while keeping uniform bounds on the geometry. This requires arguments which are more involved than the ones needed for Theorem~\ref{thm:gl connected}, and are carried out in Section~\ref{sec:isotopy bounds}.

The paper is organised as follows. In Section~\ref{sec:GL}, we define the notion of GL-sum of Riemannian manifold and prove some useful results about it.  Section~\ref{sec:cn metric surgery} is devoted to metric surgery. Theorem~\ref{thm:gl connected} is proven in Section~\ref{sec:straight}. In Section~\ref{sec:surgical_solution} we discuss surgical solutions of Ricci flow. Finally, Sections~\ref{sec:isotopy bounds} and~\ref{sec:isotopy_to_GL-can} contain the proof of Theorem~\ref{thm:isotopy_to_GL}.

\paragraph{Acknowledgments}
This research was partially supported by the Agence Nationale de la Re\-cher\-che through project~GTO~ANR-12-BS01-0014. The second author is supported by the ERC Advanced Grant 320939, GETOM (GEometry and Topology of Open Manifolds). The fourth author was partially supported by the National Science Foundation grant NSF-DMS-1509027.

\section{GL-sums} \label{sec:GL}
 
 In this section we give some details on the Gromov-Lawson connected sum construction from~\cite{Gro-Law:scalar}, (see also~\cite{marques:deforming}) which we adapt to our context of non-compact manifolds. 
 
\subsection{Topological aspects}\label{subsec:top}

Let $M$ be a (possibly disconnected) oriented 3-manifold. Let $\{(p_\alpha^-,p_\alpha^+)\}_\alpha$ be a finite or countable family of pairs of points of $M$. Assume that all these points are distinct and denote by $\mathcal{P}$ the set of all these points. Further assume that every point of $\mathcal{P}$ is isolated. Choose a family of pairwise disjoint neighbourhoods $U_\alpha^\pm$ around these points, each of which is diffeomorphic to the closed 3-ball. Denote by $M_\#(\{ p_\alpha^\pm\},\{U_\alpha^\pm\}) $ the manifold obtained by removing the interior of the $U_\alpha^\pm$'s and for each $\alpha$, gluing together $\partial U_\alpha^-$ and $\partial U_\alpha^+$ along an orientation-reversing diffeomorphism. Since the space of orientation-reversing diffeomorphisms of $S^2$ is path-connected, the diffeomorphism type of this manifold does not depend on the choice of gluing diffeomorphisms. We refrain from calling $M_\#(\{ p_\alpha^\pm\},\{U_\alpha^\pm\}) $ a `connected sum', since it needs \emph{not} be connected. 

\begin{rem}
This construction is more general than the notion of connected sum along a graph from~\cite{B2M:scalar} described in the introduction of this paper in two respects: the pieces need not be compact, and the resulting manifold need not be connected. The motivation for this greater generality comes from surgical solutions: this operation is the one needed to reconstruct the manifold before surgery from the manifold after surgery in the backward inductive arguments in Section~\ref{sec:isotopy_to_GL-can}.
\end{rem}

The manifold $M_\#(\{ p_\alpha^\pm\},\{U_\alpha^\pm\}) $  may a priori depend on the choice of the points $\{ p_\alpha^\pm\}$ and the open sets $\{U_\alpha^\pm\}$. First we show that, up to diffeomorphism, it does not depend on $\{U_\alpha^\pm\}$:
  
\begin{lem}\label{lem:independence1}
Let $M$ be an oriented 3-manifold and  $\{ p_\alpha^\pm\}$ be as above. Let $\{ U_\alpha^\pm\}$ (resp.  $\{ V_\alpha^\pm\}$) be a family of pairwise disjoint neighbourhoods of the points  $p_\alpha^\pm$, each of which is diffeomorphic to a closed $3$-ball. Then $M_\#(\{ p_\alpha^\pm\},\{U_\alpha^\pm\}) $ is diffeomorphic to $M_\#(\{ p_\alpha^\pm\},\{V_\alpha^\pm\}) $.
\end{lem}

\begin{proof}
Without loss of generality we may assume that $V_\alpha^\pm\subset U_\alpha^\pm$ for every $\alpha$. By Alexander's theorem, $U_\alpha^\pm\setminus \stackrel{\circ}{V_\alpha^\pm}$ is diffeomorphic to $S^2\times [0,1]$. Hence  $M\setminus\cup  U_\alpha^\pm$ is diffeomorphic to $M\setminus\cup  V_\alpha^\pm$.
\end{proof}

As a result, we will sometimes abuse notation by writing $M_\#(\{p_\alpha^\pm\})$ for $M_\#(\{ p_\alpha^\pm\},\{U_\alpha^\pm\})$. 

By contrast, the diffeomorphism type of $M_\#(\{ p_\alpha^\pm\})$ may in fact depend on the choice of base points. We illustrate this on the following example.

\begin{example}
Let $M$ be the disjoint union of $M_0=S^2\times\rr$ and an infinite sequence $\{X_i\}_{i\ge 1}$ of copies of $RP^3$. Choose a point $p_i^-$ in each $X_i$ and a sequence of distinct points $p_i^+$ in $M_0$ exiting every compact set. Let us denote by $L$ (resp.~$R$) the half-cylinder $S^2\times (-\infty,0]$ (resp.~$S^2\times (0,+\infty)$.) We distinguish three cases:
\begin{itemize}
\item If all  but finitely many of the $p_i^+$'s are in $L$, we obtain a manifold $M_\#^1$.
\item If all but finitely many of the $p_i^+$'s are in $R$, we obtain a manifold $M_\#^2$.
\item Otherwise we call the result $M_\#^3$.
\end{itemize}
Then $M_\#^1$ is not diffeomorphic to $M_\#^3$, since $S^2\times [0, +\infty)$ properly embeds into $M_\#^1$ but not into $M_\#^3$. The manifolds $M_\#^1$ and $M_\#^2$ are diffeomorphic, although Proposition~\ref{prop:dependence2} below does not apply.
\end{example}

We now show that the construction does not depend on the choice of basepoints under some additional hypothesis:
\begin{prop}\label{prop:dependence2}
Let $M$ be an oriented 3-manifold and $\{ p_\alpha^\pm\}$ and $\{ q_\alpha^\pm\}$ be two families of pairs of points as above. Suppose that there is a locally finite family of pairs of continuous curves $\{\gamma_\alpha^\pm \}$, parametrised by $[0,1]$ and such that for each $\alpha$, we have $\gamma_\alpha^+(0)=p_\alpha^+$, $\gamma_\alpha^+(1)=q_\alpha^+$,  $\gamma_\alpha^-(0)=p_\alpha^-$, and $\gamma_\alpha^-(1)=q_\alpha^-$.Then $M_\#(\{ p_\alpha^\pm\})$ is diffeomorphic to $M_\#(\{ q_\alpha^\pm\})$.
\end{prop}

\begin{proof}
Without loss of generality we may assume that the curves $\gamma_\alpha^\pm$ are either constant maps or embeddings. Furthermore we may assume that they are pairwise disjoint except maybe at the endpoints. By symmetry it is sufficient to prove that $M_\#(\{ p_\alpha^\pm\})$ and $M_\#(\{ \gamma_\alpha^\pm (1/2)\})$ are diffeomorphic.

Let $\{ W_\alpha^\pm\}$ be a family of closed pairwise disjoint neighbourhoods of $\gamma_\alpha^\pm ([0,1/2])$ diffeomorphic to a closed $3$-ball. It is clear then that $M_\#(\{ p_\alpha^\pm\},\{ W_\alpha^\pm\})$ is diffeomorphic to $M_\#(\{\gamma_\alpha^\pm (1/2)\},\{ W_\alpha^\pm\})$.
\end{proof}

\subsection{GL-parameters and GL-sums} \label{subsec:gl}

We present below the construction of  connected sums with a precise control on the scalar curvature. We restrict the discussion to $3$-manifolds, although part of the construction works in any dimension $n\geq 3$.

\begin{defi}
Let $(M,g)$ be a Riemannian $3$-manifold, $p$ be a point of $M$ and $\{e_k\}$ be an orthonormal basis of $T_p M$. A triple $(\rho_0,\sigma,\eta)$ of positive real numbers is called a \emph{set of GL parameters at} $p$ (with respect to~$g$) if the following requirements are met:
\begin{enumerate}
\item $\rho_0 \leq \min\{\frac{1}2 \mathrm{inj}_M(p),1\}$;
\item The scalar curvature of $g$ is greater than $\sigma$ on the ball $B_{\rho_0}(p)$;
\item The $C^2$-norm of $g$ in exponential coordinates in the basis $\{e_k\}$ is bounded from above by $\eta^{-1}$.
\end{enumerate}
\end{defi}

The motivation comes from the following result, which is essentially due to Gromov and Lawson, although the precise control of the geometry is not made explicit in their paper.

 \begin{prop}\label{prop:GL}
 For all $\rho_0,\sigma,\eta>0$ there exists $\rho_2=\rho_2(\rho_0,\sigma,\eta) \in (0,\rho_0)$ such that the following holds. 
 Let $(M,g)$ be a Riemannian $3$-manifold of positive scalar curvature, $p$ be a point of $M$ and $\{e_k\}$ be an orthonormal basis of $T_pM$. If $(\rho_0,\sigma,\eta)$ is a set of GL-parameters at $p$,  then there is a metric $g'$ on $B_{\rho_0}(p) \setminus \{p\}$ with the following properties: 
 \begin{enumerate} 	
 \item The metric $g'$ coincides with $g$ near $\partial B_{\rho_0}(p)$, 
 \item $(B_{\rho_2}(p)\setminus \{p\},g')$ is isometric to a half-cylinder,
 \item The scalar curvature of $g'$ is greater than $\frac{9}{10}\sigma$.
 \end{enumerate}
 Moreover,  if $|D^k \Rm_{g}| \leq C_k$ on $B_{\rho_0}(p)$ for $k\in\{0,1,\ldots,\bar k\}$, then $|D^k \Rm_{g'}| \leq C'_k$ on $B_{\rho_0}(p) \setminus \{p\}$, where $C'_k$ depends on $C_0,\ldots,C_k$ and $\rho_2,\rho_0,\sigma,\eta$. \\
 \end{prop}

 \begin {proof} The main step is the construction of a submanifold $M'$ of the Riemannian product $B_{\rho_0}(p) \times \RR$. This manifold is obtained by revolution along the 
 $\RR$ axis of a carefully chosen planar curve $\gamma \subset \RR^2$.  We identify  $B_{\rho_0}(p)$ with $B_{\rho_0}(0) \subset \RR^3$ using exponential normal  coordinates in the basis $\{e_k\}$. This being fixed, \cite{Gro-Law:scalar} proves the existence of a curve $\gamma \subset \RR^2$ with the following properties:
 \begin{enumerate}[(1)]
  \item the image of $\gamma$ is contained in the region $\{(\rho,t) : \rho \geq 0, t\geq 0\}$;
  \item the image of $\gamma$ contains the horizontal half-line $\rho \geq \rho_1$, $t=0$ for some $0< \rho_1 < \rho_0$;
  \item the image of $\gamma$ contains the vertical half-line $\{\rho=\rho_2, t \geq t_2\}$ for some $0<\rho_2< \rho_1$ and $t_2>0$.
  \item the induced metric on $M'=\{(x,t) :(|x|,t) \in \gamma\}$ as a submanifold of the Riemannian product $B_{\rho_0}(p) \times \RR$ has positive scalar curvature. 
 \end{enumerate}
 The curve $\gamma$ satisfying (1)-(4) is not unique. Inspecting the construction in \cite{Gro-Law:scalar} pages 425-429, we see that the choice of $\gamma$ can be refined so that $M'$ has scalar curvature greater than $\frac{9}{10}\sigma$.  
 
By choosing the radius $\rho_2$ small enough depending on $\eta$, the induced metric on $M' \cap \{t \geq t_2\}$ can be made arbitrarily close to the cylindrical metric on $S_{\rho_2}^{2}(0)\times [0, +\infty )$, where  $S_{\rho_2}^{2}(0) \subset \RR^3$ is the standard sphere of radius $\rho_2$. Using a cutoff function, the induced metric can be deformed on $M' \cap \{ t_2 \leq t \leq t_3\}$, for some  $t_3 > t_2$, into a metric $g_{M'}$ which agrees with the original metric for $t \leq t_2$ and is isometric to $S_{\rho_2}^{2}(0)\times [t_3,\infty)$ when $t \geq t_3$. Identifying $B_{\rho_0}(p) \setminus \{p\} \approx B_{\rho_0}(0)\setminus\{0\}$ and $M'$, in such a way that $S_{\rho_2}(p)=\partial B_{\rho_2}(p)$ corresponds to $M' \cap \{t=t_3+1\}$, we obtain the required metric $g'$ on $B_{\rho_0}(p) \setminus \{p\}$ by pulling back $g_{M'}$. Once $\gamma(\rho_0,\sigma,\eta)$ and the cutoff function are chosen, the estimates on the derivatives of the curvature of $g'$ follow from the definition of $M'$ as a submanifold of the product $B_{\rho_0}(p) \times \RR$.  
  \end{proof}
  
Throughout the paper, we fix for each triple $(\rho_0,\sigma,\eta)$ a curve $\gamma=\gamma(\rho_0,\sigma,\eta)$ satisfying properties (1)--(4) in the proof of Proposition~\ref{prop:GL}, in such a way that for given positive numbers  $\bar \rho$, $\bar \sigma$, $\bar \eta$, whenever $\rho_0 \geq \bar \rho$, $\sigma \geq \bar \sigma$ and $\eta \geq \bar \eta$, then $g'$ is controlled in terms of the geometry of $g$ and $\bar \rho$, $\bar \sigma$, $\bar \eta$. 

\bigskip
 
We now define the notion of GL-sum. Let $C,\sigma$ be two positive constants and $(M,g)$ be an oriented Riemannian 3-manifold with geometry bounded by $C$ and scalar curvature greater than $\sigma$. Let $\{(p_\alpha^-,p_\alpha^+)\}_\alpha$ be a finite or countable family of pairs of points of $M$. For each $\alpha$ and each $\epsilon\in\{\pm1\}$ fix a positive orthonormal basis $\{e_{\alpha,k}^\epsilon\}$ at $T_{p_\alpha^\epsilon}M$. Let $\{(\rho_\alpha,\sigma_\alpha,\eta_\alpha)\}_\alpha$ be a family of triples of positive real numbers such that for each $\alpha$, the triple $(\rho_\alpha,\sigma_\alpha,\eta_\alpha)$ is a set of GL-parameters at both $p_\alpha^+$ and $p_\alpha^-$. 
Further assume that the balls $B_{\rho_\alpha}(p_\alpha^\pm)$ are all pairwise disjoint, and that $\inf \rho_\alpha =\rho>0$ and $\inf \eta_\alpha=\eta>0$ (one has also $\inf \sigma_\alpha \geq \sigma$). For each $\alpha$, let $\rho_{2,\alpha}=\rho_{2,\alpha}(\rho_\alpha,\sigma_\alpha,\eta_\alpha)$ be given by Proposition~\ref{prop:GL} and denote by $U_\alpha^\pm$ the closure of $B(p_\alpha^\pm, \rho_{2,\alpha})$.

Then the \emph{GL-sum} associated to the above data is the Riemannian 3-manifold $(M_\#,g_\#)$, where $M_\#=M_\#(\{ p_\alpha^\pm\},\{U_\alpha^\pm\})$ is the manifold defined in subsection~\ref{subsec:top}, and $g_\#$ is as follows: for each $\alpha$, we glue together $(B_{\rho_0}(p_\alpha^+) \setminus B_{\rho_2}(p_\alpha^+),g')$ and $(B_{\rho_0}(p_\alpha^-) \setminus B_{\rho_2}(p_\alpha^-),g')$ (with the metric $g'$ given in each case by Proposition~\ref{prop:GL})
by identifying their boundaries, which are isometric to $S_{\rho_2}(0) \subset \RR^n$, by the orientation reversing isometry obtained by composing exponential maps and the identification of the tangent spaces given by the linear map sending $\{ e_{\alpha,1}^+,e_{\alpha,2}^+,,e_{\alpha,3}^+\}$ to $\{ - e_{\alpha,1}^-,e_{\alpha,2}^-,,e_{\alpha,3}^-\}$. 

Sometimes we will use the notation $(M,g)_\#$, for $(M_\#,g_\#)$. When we need to be more precise, we will specify the basepoints and use the notation $(M_\#(\{p_\alpha^\pm\}),g_\#(\{p_\alpha^\pm\})$ or $(M,g)_\#(\{p_\alpha^\pm\})$, or even indicate all parameters:  $(M,g)_\#(\{(\rho_\alpha,\sigma_\alpha,\eta_\alpha)\},\{p_\alpha^\pm\},\{e_{k,\alpha}^\pm\})$.  In the special case of a classical connected sum $M_1\# M_2$ we use the notation $g_1\# g_2$.

\begin{rem}
By construction, each sphere along which some gluing has been done is straight with respect to $g_\#$. As a consequence, if $(M_\#,g_\#)$ is a GL-sum of closed manifolds, then $g_\#$ is a GL-metric (cf.~Definition~\ref{defi:GL metric}.)
\end{rem}

Observe that the metric $g_\#$ has scalar curvature $\geq \frac{9}{10} \sigma$, and geometry bounded by a constant $C'$ depending only on $C, \rho, \sigma, \eta$. However, the bounded geometry of $(M,g)$ alone does not suffice to control $\eta$. In order to do this, we need bounds on the derivative of the curvature. 

From Hamilton \cite[Corollaries~4.11 and~4.12]{hamilton:compactness} we have the following: 
\begin{lem}\label{lem:eta} For each $n\ge 2$ there is a universal constant $c(n)\in (0,1/2)$ with the following property. Let $C_0,C_1,C_2$ be positive numbers and let $(M,g)$ be an $n$-dimensional Riemannian manifold such that $|D^k \Rm| \leq C_k$ for $k=0,1,2$. Then the $C^2$-norm of $g$ in exponential coordinates is bounded from above in  $B_{c(n)C_0^{-1/2}}(p)$ by a constant depending only on $C_0,C_1,C_2$.
\end{lem}

Therefore when $|D^k \Rm| \leq C_k$, $k=0,1,2$, and $\rho_\alpha <c(3)C_0^{-1/2}$, we have a positive lower bound for $\eta$.
In the context of a Ricci flow with bounded curvature, the extra bounds on the derivatives of the curvature are provided by the so-called Shi estimates~\cite{shi:complete3}.  We note that the GL-sum construction preserves bounds on the derivatives of the curvature, as remarked in Proposition~\ref{prop:GL}.  Therefore starting from a metric with bounds on the derivatives of the curvature, Lemma~\ref{lem:eta} and Proposition~\ref{prop:GL} imply that one can iterate the GL construction keeping the geometry under control.

\subsection{Continuity of the  GL-sum construction}

We give below a generalisation of Proposition 6.1. of~\cite{marques:deforming} that we will need.

Let $(M,g)$ be a  Riemannian  3-manifold and $\bar k$ be a natural number. Given a finite family $A_k$, $k \in \{0,1,\ldots,\bar k\}$ of positive real numbers, we denote $\cA_{\bar k}=(A_k)_{k \in \{0,1,\ldots,\bar k\}}$ and we say that $g$ has geometry bounded by $\cA_{\bar k}$ if $\inj (g) \geq A_0^{-1/2}$ and $\vert D^k \Rm(g) \vert \leq A_k$ for $k \in \{0,1,\ldots,\bar k\}$. In general we will take $\bar k=[\epsilon^{-1}]$ and denote 
by $\bar \cA=\cA_{[\epsilon^{-1}]}$, omitting the index. Given $n \leq \bar k$ and  $\cA_{\bar k}$, then $\cA_n:=(A_k)_{k \in \{0,1,\ldots,n\}}$. With these notations, when $g$ has geometry bounded by $\cA_2$, Lemma \ref{lem:eta} gives a lower bound $\eta(\cA_2)>0$ for the parameter $\eta$ in balls of radius less than  $c(3){A_0}^{-1/2}$.

\begin{prop}\label{prop:continuousGL} Let $\bar k$ be a natural number. For every $\cA_{\bar k}=(A_k)_{k \in \{0,1,\ldots,\bar k\}}$,  $\rho \in (0, \min\{c(3){A_0}^{-1/2},1\})$, $\sigma>0$, and $\eta>0$, there exist $\cA^{\#}_{\bar k}=\cA^{\#}_{\bar k}(\cA_p,\rho,\sigma,\eta)=(A^{\#}_k)_{k \in \{0,1,\ldots,\bar k\}}$ with the following property. Let $M$ be an oriented  $3$-manifold and $g_t$ be a continuous path of metrics with scalar curvature $\geq \sigma$ and geometry bounded by $\cA_{\bar k}$. Let $p_{\alpha,t}^\pm$ be continuous paths of points of $M$, and $\{e_{k,\alpha}^\pm(t)\}$ be continuous paths of positive orthonormal bases of $(T_{p_{\alpha,t}^\pm}M,g_{t})$. Let $(\rho_\alpha, \sigma_\alpha,\eta_\alpha)$ be a triple which is for each $t$ a set of GL-parameters at $p_{\alpha,t}^\pm$ w.r.t.~$g_t$. Suppose that $\rho_\alpha \geq \rho,$, $\sigma_\alpha \geq \sigma$, $\eta_\alpha \geq \eta$ and that there exist pairwise disjoint closed $3$-balls $\{W_{\alpha}^\pm\}$ such that $B_{g_t}(p_{\alpha,t}^\pm,2\rho_\alpha) \subset W_\alpha^\pm$  for all $t$. Let $(M_{\#,t},g_{\#,t})$ denote the GL-sum $(M,g_t)_\#( \{(\rho_\alpha,\sigma_\alpha,\eta_\alpha)\}, \{p_{\alpha,t}^\pm\},\{e_{k,\alpha}^\pm(t)\})$. Then there exist positive diffeomorphisms $\phi_t : M_{\#,0} \to M_{\#,t}$ such that
\begin{enumerate}[\hspace{0.3cm} 1)]
\item $\phi_t$ is the identity on $M \setminus \cup W_{\alpha}^\pm \subset M_{\#,0} \cap M_{\#,t}$.
\item The pulled-back metrics $\phi_t^\ast(g_{\#,t})$ define a continuous path  on $M_{\#,0}$, with scalar curvature $\geq \frac{9}{10}\sigma$ and geometry bounded by $\cA^{\#}_{\bar k}$.
\item We have $\phi_t^\ast(g_{\#,t})=g_t$ on  $M \setminus \cup W_{\alpha}^\pm \subset M_{\#,0}$.
\end{enumerate}
  \end{prop}

 \begin{proof}
It is enough to argue in a pair $W_\alpha^{^\pm}$. By standard arguments there is an isotopy $h_t$ from $ W_\alpha^\pm$ into itself which is the identity on $\partial W_\alpha^\pm$, and such that 
\begin{enumerate}[\hspace{0.3cm}  $\bullet$]
\item $h_t$ maps $p_{\alpha,0}^\pm$ to $ p_{\alpha,t}^\pm$, 
\item $(h_t)_\ast$ sends $\{e_{k,\alpha}^\pm(0)\}$ to $\{e_{k,\alpha}^\pm(t)\}$, 
\item $h_t$ coincides with $\exp_t \circ (h_t)_\ast \circ \exp_0^{-1}$ from $B_{g_0}(p_{\alpha,0}^\pm,\rho_\alpha)$ 
to $B_{g_t}(p_{\alpha,t}^\pm,\rho_\alpha)$. 
\end{enumerate}
The isotopy is compatible with the identifications of the GL-sum construction and yields the required family of diffeomorphisms.
 \end{proof} 
 
 \begin{rem}
 It follows that $t \mapsto g_{\#,t}$ is a continuous path in $\cR_{9\sigma/10}^{\cA^{\#}_{\bar k}}(M_{\#,0})/\mathrm{Diff}(M_{\#,0})$. Another consequence is that up to isotopy and diffeomorphism, a GL-sum $(M,g)_\#(\{(\rho_\alpha,\sigma_\alpha,\eta_\alpha)\},\{p_\alpha^\pm\},\{e_{k,\alpha}^\pm\})$ does not depend on the positive orthonormal basis $\{e_{k,\alpha}^\pm\}$. Therefore in the sequel we will frequently drop the mention to the basis, which implicitely will always be positive and orthonormal. We observe that taking 
 a negative basis  may change the diffeomorphism type of the resulting manifold.
 \end{rem}
   
\subsection{Two lemmas}
For later reference we collect two technical results from~\cite{marques:deforming}.

 \begin{lem}[{cf.~\cite[p.~835]{marques:deforming}}] \label{lem:conf-gamma} Let $B^-$ and $B^+$ be balls of constant sectional curvature in $(0,A]$, and let $p^\pm \in B^\pm$. Let $\gamma_1=\gamma(\rho_1,\sigma_1,\eta_1)$ and $\gamma_2=\gamma(\rho_2,\sigma_2,\eta_2)$ be two planar curves, where $(\rho_i,\sigma_i,\eta_i)$ are GL-parameters suitable for both balls, $i=1,2$. Then the GL-sums $B^-\#_{\gamma_1} B^+$ and  
 $B^- \#_{\gamma_2} B^+$ performed at $p^\pm$ are isotopic, without changing the metric near $\partial B^{\pm}$, through metrics of scalar curvature $\geq \bar \sigma = \min\{\frac{9}{10}\sigma_1,\frac{9}{10}\sigma_2\}>0$ and geometry bounded by $\bar \cA=\bar \cA(A,\gamma_1,\gamma_2)<\infty$. 
 \end{lem}
 
The following lemma shows that, applying \cite[Prop.~3.3]{marques:deforming}, we can continuously deform a standard cylinder in such a way that in its middle it becomes isometric to a subset of a round sphere.  

 \begin{lem}[{cf.~\cite[Prop.~3.3]{marques:deforming}}]\label{lem:conf-cyl} Let $(S^3,h)$ be a standard $3$-sphere of scalar curvature $\sigma$ and let $x \in S^3$. Denote by $h'$ the metric obtained on $S^3 \setminus \{-x,x\}$ by applying the GL-sum construction with parameter 
 $(\rho,\sigma,\eta)$ at both $x$ and $-x$, and let $\bar \sigma'$ be the (constant) scalar curvature of $h'$ near $\pm x$. Let $(S^2 \times (a,b),g)$ be a standard cylinder of scalar curvature $\bar \sigma'$, and let 
 $(a',b') \subset [a,b]$. Then there exists an isotopy $(g_t)$ on the cylinder such that 
 \begin{enumerate}
 \item $g_t=g$ on $S^2 \times (a,a')$ and $S^2 \times (b',b)$.
 \item $g_0=g$ and $(S^2 \times (a',b'),g_1)$ is isometric to a subset of $(S^3 \setminus \{-x,x\},h')$.
 \item Each metric $g_t$ has scalar curvature greater than $\bar{\bar \sigma} = \bar{\bar \sigma}(\rho,\sigma,\eta)>0$ and geometry bounded by $\bar{\bar \cA}= \bar{\bar \cA}(\rho,\sigma,\eta)<\infty$.
 \end{enumerate}
  \end{lem}

\section{Metric surgery}\label{sec:cn metric surgery}
In this section, we recall some notions and results from~\cite{B2M:scalar}.
Throughout the paper we denote by $d\theta^2$ the round metric of scalar curvature $1$ on $S^2$, and whenever $I\subset\RR$ is an open interval, we denote by $\gcyl$ the product metric $d\theta^2 + ds^2$ on $S^2\times I$. We also use this notation on $\RR^3\setminus \{0\}$, working in polar coordinates. The origin $0$ of $\RR^3$  is called the \bydef{tip} and the radial coordinate is denoted by $r$.

A \emph{standard initial metric} is a complete metric on $\RR^3$ which is 
rotationally symmetric, has bounded nonnegative sectional curvature, is isometric to $\gcyl$ on some neighbourhood of infinity and is round on some neighbourhood of the tip. As in~\cite[section 3.2]{B2M:scalar} (cf.~also~\cite[section 7.1]{B3MP}), we fix a particular standard initial metric $g_u$ which is cylindrical outside $B(0,3)$ and round of sectional curvature $1/2$ near the tip. We also fix a smooth, nonincreasing function $f:[0,+\infty)\to[0,+\infty)$ with support in $[0,5]$, and set $g_{\std}  = e^{-2f(r)}g_u$. The function $f$ is chosen so that metric surgery is distance-nonincreasing and preserves the so-called  Hamilton-Ivey pinching property, none of which is important in the present paper. The metric $g_\std$ is also a standard initial metric. It has scalar curvature $\geq 1$ everywhere and positive sectional curvature on $B(0,5)$.

Let $g$ be a Riemannian metric on  $S^2 \times (-4,4)$ which is $\epsi$-close to the metric $\gcyl$  in the $C^{[\epsi^{-1}]}$-topology. We now describe a surgery operation which turns $(S^2 \times (-4,4),g)$ into a Riemannian manifold $(\cS^- \sqcup \cS^+,g_\surg)$, where $\cS^- \sqcup \cS^+$ is a disjoint union of copies of the open ball $B(0,9) \subset \RR^3$ and the metric $g_\surg$ is defined as follows. Let $\psi^- : S^2\times (-4,4) \to B(0,9)$ be the embedding given by $(\theta,s) \mapsto ((5-s),\theta)$ in polar coordinates. Fix a function $\chi : \bar B(0,9) \to [0,1]$ such that $\chi \equiv 0$ on $B(0,3)$ and $\chi \equiv 1$ outside $B(0,4)$.


Let $g_\surg^-$ be the Riemannian metric on $\cS^-$ defined as follows:
$$ g_\surg^- = \left\{ \begin{array}{ll}
g_\std & \text{on } B(0,3) \\
\chi e^{-2f}(\psi^-)_\ast g  + (1-\chi) g_\std & \text{when } 3 \leq r \leq 4\\
e^{-2f}(\psi^-)_\ast g & \text{when } 4 \leq r \leq 5 \\
(\psi^-)_\ast g   & \text{when } 5 \leq r \leq 9.
\end{array}
\right.
$$

In the case that $g=\gcyl$, 
one can check that $g_\surg^- =g_\std$. More generally,  $g_\surg^-$ is $\delta'(\epsi)$-close to $g_{\std}$ for some $\delta'(\epsi)$ going to zero with $\epsi$.  Notice that the metric on $S^2\times (-4,0)$ remains unchanged (up to a diffeomorphism). This construction thus amounts to capping off an `almost standard cylinder'.

Likewise, we define a metric $g_\surg^+$ on $\cS^+$ using the embedding  $\psi^+ : S^2\times (-4,4) \to B(0,9)$ given in polar coordinates by $(\theta,s) \mapsto ((5+s),\theta)$. Finally we let $g_\surg$ be the metric on $\cS^- \sqcup \cS^+$ whose restriction to $\cS^-$ (resp.~to~$\cS^+$) is $g_\surg^-$ (resp.~$g_\surg^+$.)

From~\cite[section 3.3]{B2M:scalar}, we get the following result. Let $\epsi_0>0$ be the number defined in Lemma~3.6 of~\cite{B2M:scalar}.

\begin{theo}\label{thm:surgery} There exists $\delta_0 \in (0,\frac{\epsi_0}{10})$ and a function $\delta': (0,\delta_0] \to (0,\epsi_0/10)$ going to zero at zero with the following property. Let $\delta \in (0,\delta_0)$ and let $g$ be a metric on $S^2 \times (-4,4)$ which is $\delta$-close to $\gcyl$. Then the Riemannian manifold $(\cS^-,g_\surg^-)$ (resp.~$(\cS^+,g_\surg^+)$) has the following properties:
	\begin{enumerate}
		\item All sectional curvature are positive on $B(0,4)$.
		\item The scalar curvature is $\geq R_g$ on $r^{-1}((4,9))$.
		\item The smallest eigenvalue of the curvature operator of $g_\surg^\pm$ is larger than or equal to the smallest eigenvalue of $\Rm_g$ on $r^{-1}((4,9))$. 		
		\item The metric $g_\surg^\pm$ is $\delta'(\delta)$-close to $g_{\std}$.
	\end{enumerate}
\end{theo}

Our next goal is to define metric surgery on a neck in a Riemannian manifold. For technical reasons, it is useful to allow the length of the neck to vary. 

\begin{defi}\label{defi:neck}
Let $\epsi,L$ be positive real numbers and $(M,g)$ be a Riemannian 3-manifold. An $(\epsi , L)$-\emph{neck} in $M$ is an open subset $N\subset M$ for which there is a $C^{[\epsi^{-1}]+2}$-diffeomorphism $\phi : S^2 \times (-L,L) \to N$, called a \bydef{parametrisation,} and a number $\lambda>0$ such that $\lambda \phi^\ast g$ is $\epsi$-close to $\gcyl$ in the $C^{[\epsi^{-1}]+2}$-topology (defined by $\gcyl$.) The set $\phi(S^2 \times \{0\})$ is called the \emph{middle sphere} of $N$. When $L=\epsi^{-1}$ we simply say that $N$ is an $\epsi$-\emph{neck}.
\end{defi}

We recall that $\epsi$-closeness in the $C^{[\epsi^{-1}]+2}$-topology between the metrics yields a control on the derivatives of the curvature up to order $[\epsi^{-1}]$. In the sequel, $\epsi$-closeness will always be understood in this topology.

Let $(M,g)$ be a Riemannian $3$-manifold, and $N\subset M$ be an $(\epsi, 4)$-neck. Let $\phi, \lambda$ be as above and $S$ be the middle sphere of the neck. We call \emph{metric surgery on $N$} (or along $S$) the procedure of replacing $(M,g)$ by the Riemannian manifold $(M_+,g_+)$, where
$$ M_+ = \left((M\,\setminus S) \,\, \sqcup \,\, \cS^- \sqcup \cS^+\right)/\sim $$
identifying $\phi(S^2 \times (-4,0])$
with $\cS^- \setminus B(0,5)$ via $\psi^- \circ \phi^{-1}$ and $\phi(S^2 \times [0,4))$ with $\cS^+ \setminus B(0,5)$ via $\psi^+ \circ \phi^{-1}$, and 
$$ \left\{\begin{array}{ll}
g_+=g & \textrm{ on } M \setminus S\\
g_+=\lambda^{-1}(\lambda \phi^\ast g)_\surg & \textrm{ on }  \cS^- \sqcup \cS^+.
\end{array} \right.
$$

\begin{rem}
If $g$ is $C^{\infty}$-smooth, then the manifold $(M_+, g_+)$ is  $C^{[\epsi^{-1}]+2}$-smooth. Furthermore if $g$ satisfies $\vert D^kRm\vert\leq A_k$ for  $k\leq [\epsi^{-1}]$,  then $g_+$ satisfies $\vert  D^kRm\vert\leq B_k$ for $k\leq [\epsi^{-1}]$, where $\cB = \{B_k\}$ depends on $\cA= \{A_k\}$ only.
\end{rem}

\begin{rem}
The value $4$ above can be replaced by any positive value $L$ (replacing $\cS^{\pm}$ with $B(0,5+L)$ if $L<4$). In particular this works for $L=\epsi^{-1}$, so that we have  defined the notion of surgery on an $\epsi$-neck .
\end{rem}

In some sense, the metric surgery process can be reversed by the GL-sum construction, as shown by the following result. Choose orientations on the manifolds $S^2 \times (-4,4)$, $\cS^-$ and $\cS^+$ in such a way that the diffeomorphisms $\psi^-$ and $\psi^+$ are positive. Let $\cS^-\#\cS^+$ be the connected sum with basepoints the tips. We choose once and for all an identification between $S^2\times (-4,4)$ and $\cS^-\#\cS^+$ which coincides with $\psi^-$ and $\psi^+$ near the boundary (as in the proof of~\cite[Lemma~6.2]{marques:deforming}.) This does not depend on the choice of orientations. Let $\epsi_3$ be the
constant given by~\cite[Lemma~6.3]{marques:deforming}. From this lemma we get:
\begin{lem}\label{lem:gl-surg}
Let $\epsi\in (0,\epsi_3)$. Let $g$ be a Riemannian metric on $S^2 \times (-4,4)$ which is $\epsi$-close to $\gcyl$ on $S^2 \times (-4,4)$. Then $(\cS^- \# \cS^+, g_{\surg}^-\#g_{\surg }^+)$ can be continuously deformed back
	into $(S^2 \times  (-4,4), g)$ through metrics which all coincide
	with $g$ near the ends of $S^2 \times  (-4,4)$, have scalar curvature greater than $9/10$ and geometry bounded by some $\cB_{[\epsi^{-1}]}$. 	
\end{lem} 	
If $N$ is an $(\epsi , 4)$-neck we use its parametrisation to identify it with $\cS^-\#\cS^+$.
Finally if $(M_+,g_+)$ is obtained from $(M,g)$ by surgery in a family of disjoints necks, then we can identify $(M_+)_\#$ with $M$ canonically. This identification is the identity in the complement of the necks. The above lemma shows that $(g_+)_\#$ is isotopic to $g$.

\section{Connecting two GL-metrics: proof of Theorem \ref{thm:gl connected}} \label{sec:straight} \label{sec:mainproofs}
The aim of this section is to prove Theorem \ref{thm:gl connected}. This theorem follows immediately from  Lemma~\ref{lem:same cS isotopic} and~Proposition~\ref{prop:isotopy-straight2} below, which are proven in Subsections~\ref{subsec:same cS isotopic}
and \ref{subsec:isotopy-straight2} respectively.

Given a spherical splitting $\cS=\{S_\alpha\}$ of a $3$-manifold $M$ and two metrics $g,h$ on $M$, we will say that $g=h$ \emph{near} $\cS$ if for every $\alpha$ there is a neighbourhood of $S_\alpha$ on which $g$ and $h$ coincide. 
\begin{lem} \label{lem:same cS isotopic} Let $M$ be a non-compact $3$-manifold. Let $(g,h,\cS)$ be a triple such that $g$ and $h$ are GL-metrics on $M$ that belong to $\cR_1^{bg}(M)$, and $\cS$ is a spherical splitting which is straight for both $g$ and $h$. Suppose that $g=h$ near $\cS$. Then $g$ is  isotopic to $h$ modulo diffeomorphism.
\end{lem}

\begin{prop}\label{prop:isotopy-straight2} Let $M$ be a non-compact $3$-manifold, let $g,g' \in \cR_1^{bg}(M)$ be two GL-metrics on $M$ and let $\cS'$ be a spherical splitting straight for $g'$. Then there exists a metric $h \in \cR_1^{bg}(M)$, isotopic to $g$ modulo diffeomorphism, such that $h=g'$ near $\cS'$. (In particular, $\cS'$ is straight for $h$.)
\end{prop}

\subsection{Proof of Lemma~\ref{lem:same cS isotopic}}
\label{subsec:same cS isotopic}
We start with a lemma which allows to combine paths of metrics defined on compact submanifolds in an exhaustion. In the sequel, the \bydef{support} of a path of metrics $g_t$ on a manifold $M$ is the closure of the subset of $M$ where $t\mapsto g_t$ fails to be constant. 

\begin{lem}\label{lem:combining paths}
	Let $M$ be a 3-manifold and $C$ be a positive constant. Let  $\emptyset =K_{-1}\subset K_0\subset K_1\subset K_2\subset \cdots K_\ell \cdots$ be an exhaustion of $M$ by compact submanifolds. Suppose we are given a sequence of continuous paths $h^{(\ell)}:[\frac{\ell}{\ell+1},\frac{\ell+1}{\ell+2}]\to \mathcal{R}_{1}^{bg}(M)$ such that
	\begin{enumerate}[(i)]
		\item For each $\ell$ we have $h^{(\ell)}(\frac{\ell+1}{\ell+2})=h^{(\ell+1)}(\frac{\ell+1}{\ell+2})$;
		\item For each $\ell$, the path $h^{(\ell)}$ has support in $K_{\ell}\setminus K_{\ell-1}$;
		\item For each $\ell$, the metric $h^{(\ell)}(\frac{\ell+1}{\ell+2})$ has geometry bounded by $C$.
	\end{enumerate}
	Then there is a unique continuous path $h:[0,1] \to \mathcal{R}_{1}^{bg}(M)$ which extends the $h^{(\ell)}$'s. Furthermore, $h(1)$ has geometry bounded by $C$.
\end{lem}

\begin{proof}
	By~(i) there is a unique extension $h$ of the paths $h^{(\ell)}$ to the interval $[0,1)$. By~(ii), $h(t)$ has a limit in the $\mathcal{C}^{\infty}_{loc}$ topology when $t$ goes to $1$. We set $h(1)$ equal to this limit, thus obtaining a continuous map defined on $[0,1]$. Condition~(iii) ensures that $h(1)$ has geometry bounded by $C$.
\end{proof}

\begin{rem}\label{rem:combining paths}
	A typical application of Lemma~\ref{lem:combining paths} is as follows: $M$ is a disjoint union of closed manifolds $M_\ell$ for $\ell\in\nn$; we have a metric $g\in \cR_1^{bg}$, which we wish to deform to another metric with nice properties, and we are given a collection of paths $g^{(\ell)}:[0,1]\to \cR_{1}(M_\ell)$ such that $g^{(\ell)}(0)\equiv g$ on $M_\ell$ for all $\ell$, and such that all final metrics $g^{(\ell)}(1)$ have geometry uniformly bounded by some constant $C$. We first reparametrize those paths such that the domain of $g^{(\ell)}$ is $[\frac{\ell}{\ell+1},\frac{\ell+1}{\ell+2}]$ for every $\ell$. We then extend them  to continuous paths $h^{(\ell)}:[\frac{\ell}{\ell+1},\frac{\ell+1}{\ell+2}]\to \cR_1^{bg}(M)$ such that each $h^{(\ell)}$ is constant outside $M_\ell$. The fact that we stay in $\cR_1^{bg}(M)$ is due to the hypothesis that $g$ has bounded geometry together with the fact that each path has compact support. Note that we are \emph{not} claiming that the paths have image in $\cR_1^{C'}$ for a uniform $C'$. We then set $K_\ell:=\bigcup_{k\le \ell} M_k$ for all $\ell$. We then apply Lemma~\ref{lem:combining paths} in order to obtain the global deformation of $g$.
	
Note that the `obvious' way of combining the paths, without reparametrizing, could lead to a path of metrics which does not stay in $\cR_1^{bg}(M)$.
\end{rem}

In the sequel, we sometimes have paths of metrics which lie in $\cR_{9/10}^{bg}(M)$ rather than $\cR_1^{bg}(M)$, because we take GL-sums. This is not a problem because of the following lemma.

\begin{lem}\label{lem:sigma}
Let $M$ be a 3-manifold. Let $g,g'$ be two metrics in $\cR_1^{bg}(M)$. If there exists $\sigma>0$ such that $g,g'$ are isotopic in $\cR_\sigma^{bg}(M)$, then $g,g'$ are isotopic in $\cR_1^{bg}(M)$. 
\end{lem}

\begin{proof}
Let $g_t$ be an isotopy from $g$ to $g'$ in $\cR_\sigma^{bg}(M)$. For each $t$ we have $\sqrt{\sigma} g_t \in \cR_1^{bg}(M)$. Hence we can isotope $g$ linearly to $\sqrt{\sigma} g$, follow $\sqrt{\sigma} g_t$, and finally isotope $\sqrt{\sigma} g'$ to $g'$.
\end{proof}

The following technical result will allow us to obtain an isotopy on a GL-sum manifold from isotopies on its components.  
\begin{lem}\label{lem:som gl gb}
	Let $M$ be an oriented $3$-manifold which is a disjoint union of closed manifolds. Let $\{p_\alpha^\pm\}$ and $\{q_\alpha^\pm\}$ be two families of pairs of points of $M$ such that for every $\epsilon \in \{\pm\}$ and every $\alpha$, the points $p_\alpha^\epsilon$ and $q_\alpha^\epsilon$ lie in the same connected component of $M$. Let $g_t $ be a path of metrics in $\cR_1^{bg}(M)$. Let $\rho$ and $\eta$ be positive numbers such that $(\rho,1, \eta)$ is a set of GL-parameters at $\{p_\alpha^\pm\}$ w.r.t.~$g_0$ (resp. at $\{q_\alpha^\pm\}$ w.r.t.~$g_1$) and such that 
	the $2\rho$-balls for $g_0$ centred at  $\{p_\alpha^\pm\}$ (resp. the $2\rho$-balls for $g_1$ centred at  $\{q_\alpha^\pm\}$) are round and pairwise disjoint. 
	
	Then there is a continuous path of metrics $g_\#(t)$ on $M_\#=M_\#(\{(\rho,1,\eta)\},\{p_\alpha^\pm\})$ with  scalar curvature greater than $9/10$ and bounded geometry, such that $g_\#(0)= (g_0)_\#(\{(\rho,1,\eta)\},\{p_\alpha^\pm\})$ and $(M_\#, g_\#(1))$ is isometric to $(M,g_1)_\#(\{(\rho,1,\eta)\},\{q_\alpha^\pm\}))$.
\end{lem}

\begin{proof}
In the compact case this follows from Proposition \ref{prop:continuousGL}, so we assume that  $M$ is non-compact and denote by $\{M_\ell\}_{\ell \in \nn}$ its components, which by hypothesis are closed. We begin by some reductions.
	
	Up to slightly moving points of $\{p_\alpha^\pm\}$ and applying Proposition \ref{prop:continuousGL} to get a corresponding deformation of $(g_0)_\#(\{p_\alpha^\pm\})$, we can assume that for all $x \in \{p_\alpha^\pm\}$ and $y \in \{q_\alpha^\pm\}$ one has $x \not=y$.  
	
	Without loss of generality we may assume that the support of the restriction of $g_t$ to the interval 
	$[\frac{\ell}{\ell+1},\frac{\ell+1}{\ell+2}]$ is contained in $M_\ell$. As $g_0$ and $g_1$ have geometry bounded by a constant $C$, it follows that $g_t$ has geometry bounded by a constant depending on $\ell$ when 
	$t \in [\frac{\ell}{\ell+1},\frac{\ell+1}{\ell+2}]$. 
	
	We then deform $(M_\#(\{p_\alpha^\pm\}),(g_0)_\#(\{p_\alpha^\pm\}))=(M,g_0)_\#$, up to diffeomorphism, as follows. We apply Lemma \ref{lem:conf-cyl}  to deform each GL-neck $B_{g_0}(p^-,\rho)\, \# \, B_{g_0}(p^+,\rho)$, where $p^\pm \in \{p_\alpha^\pm\}$, into the GL-sum
	$$B_{g_0}(p^-,\rho)\,\, \# \,\,S^3(p^\pm) \,\, \# \,\, B_{g_0}(p^+,\rho)$$ where 
	$S^3(p^\pm)$ is a round $3$-sphere of scalar curvature $1$, and the GL-sum is made 
	at $\{(p^-,x(p^-)),(p^+,x(p^+))\}$ where $x(p^-)=-x(p^+) \in S^3(p^\pm)$. 
	This has the effect of replacing $(M,g_0)_\#$ by the GL-sum $(M \sqcup X,g_0 \sqcup g_X)_\#(\{(p_\alpha^-,x(p_\alpha^-)),(p_\alpha^+,x(p_\alpha^+))\}_\alpha)$, where $(X,g_X)$ is the disjoint union of the round $3$-spheres $S^3(p_\alpha^\pm)$. This allows to consider the GL-sums made  on a component $(M_\ell,g_t)$  independently of the other components of $M$, the component $M_\ell$ being connected to a fixed disjoint union of $3$-spheres along the path.

	We now begin the construction of the isotopy, with initial point the manifold $(M \sqcup X,g_0 \sqcup g_X)_\#$. 
	In every component $M_\ell$, choose  pairwise disjoint paths $[0,1] \mapsto p_{\alpha,t}^\epsilon$, constant outside $[\frac{\ell}{\ell+1},\frac{\ell+1}{\ell+2}]$, connecting $p_{\alpha}^\epsilon$ to $q_{\alpha}^\epsilon$ (the paths are disjoint thanks to the first simplification above; we also choose paths of positive orthonormal bases which we do not mention). 
	By compactness of $M_\ell \times [\frac{\ell}{\ell+1},\frac{\ell+1}{\ell+2}]$,  there exist $\rho^{(\ell)}>0$ and $\eta^{(\ell)}>0$ such that 
	$(\rho^{(\ell)},1,\eta^{(\ell)})$ is a set of GL-parameters at $p_{\alpha,t}^\epsilon$ with respect to every metric $g_t$ such that $t \in [\frac{\ell}{\ell+1},\frac{\ell+1}{\ell+2}]$. Moreover,  we choose $\rho^{(\ell)} $ small enough so that there exists a family of pairwise disjoint $3$-balls $W_\alpha^{\epsilon,(\ell)} \subset M_\ell$ such that 
	$B_{g_t}(p_{\alpha,t}^\epsilon, 2\rho^\ell) \subset W_\alpha^{\epsilon,(\ell)}$ for all 
	$t \in [\frac{\ell}{\ell+1},\frac{\ell+1}{\ell+2}]$. 
	
	For each  $t \in [\frac{\ell}{\ell+1},\frac{\ell+1}{\ell+2}]$, define the Riemannian manifold $(M_{\#,t}^{(\ell)},g_{\#,t}^{(\ell)})$ as the GL-sum $(M \sqcup X, g_t \sqcup g_X)_\#(\{(p_{\alpha,t}^-,x(p_\alpha^-),(p_{\alpha,t}^+,x(p_\alpha^+)\}_\alpha)$, using  GL-parameters $(\rho^{(\ell)},1,\eta^{(\ell)})$ when $p_\alpha^\epsilon \in M_\ell$ and $(\rho,1,\eta)$ when $p_\alpha^\epsilon \notin M_\ell$. From Proposition  \ref{prop:continuousGL}, each metric $g_{\#,t}^{(\ell)}$ has scalar curvature greater then $9/10$ and geometry bounded by a constant depending on $\ell$; moreover, there exist  diffeomorphisms  $\phi_{t}^{(\ell)} : M_{\#,\frac{\ell}{\ell+1}}^{(\ell)} \to M_{\#,t}^{(\ell)}$ such that the pulled-back metrics $(\phi_{t}^{(\ell)})^\ast g_{\#,t}^{(\ell)}$ define a continuous path on $M_{\#,\frac{\ell}{\ell+1}}^{(\ell)}$. The manifolds 	$M_{\#,\frac{\ell+1}{\ell+2}}^{(\ell)}$ and $M_{\#,\frac{\ell+1}{\ell+2}}^{(\ell+1)}$ differ only on a union of GL-necks connected to $M_\ell$ and to $M_{\ell+1}$, due to the difference in the GL-parameters used to define these necks: from $(\rho^{(\ell)},1,\eta^{(\ell)})$ to  $(\rho,1,\eta)$ in the case of $M_\ell$ and from  $(\rho,1,\eta)$ to $(\rho^{(\ell+1)},1,\eta^{(\ell+1)})$ in the case of $M_{\ell+1}$. We identify them by identifying the corresponding necks. Pulling back the paths $g_{\#,t}^{(\ell)}$
	by diffeomorphisms $\phi_{\frac{\ell}{\ell+1}}^{(\ell-1)}\circ \cdots \circ \phi_{\frac{1}{2}}^{(0)}$ and concatenating, we obtain on the manifold $M_{\#,0}^{(0)}$ a piecewise continuous path of metrics defined on $[0,1)$. Now, the GL-sums are made with different parameters creating discontinuities of the path at times $\frac{\ell+1}{\ell+2}$; thanks to Lemma \ref{lem:conf-gamma} they  can be smoothed out and we obtain a continuous path $g_{\#}(t)$ defined on $[0,1)$, with scalar curvature greater than $9/10$ and bounded geometry. From the construction above, for any exhaustion  
	$K_0 \subset K_1 \subset \ldots $ of  $M_{\#,0}^{(0)}$ by compact 
	subsets, $t \mapsto g_{\#}(t)$ is constant  on 
	$K_j$ for all $t$ close enough to $1$. It follows that $t \mapsto g_{\#}(t)$ extends continuously to the interval $[0,1]$. Moreover, denoting 
	$(M_{\#,1}^{(\infty)},g_{\#,1}^{(\infty)})$ the GL-sum $(M \sqcup X, g_1 \sqcup g_X)_\#(\{(\rho,1,\eta)\}, \{(p_{\alpha,1}^-,x(p_\alpha^-),(p_{\alpha,1}^+,x(p_\alpha^+)\}_\alpha)$, then the sequence of diffeomorphisms $\phi_{\frac{\ell}{\ell+1}}^{(\ell)} \circ \ldots \circ \phi_{1/2}^{(0)}$ converges to an isometry from $(M_{\#,0}^{(0)},g_{\#}(1))$  to  $(M_{\#,1}^{(\infty)},g_{\#,1}^{(\infty)})$. It is easy to further deform
	$(M_{\#,1}^{(\infty)},g_{\#,1}^{(\infty)})$ into  $(M_\#(\{q_\alpha^\pm\}),(g_1)_\#(\{q_\alpha^\pm\}))$, using Lemma 
	\ref{lem:conf-cyl}. This concludes the proof of Lemma~\ref{lem:som gl gb}.
\end{proof}

\begin{proof}[Proof of Lemma~\ref{lem:same cS isotopic}]
Fix an orientation of $M$. Let $(\widehat M, \widehat g)$ (resp. $(\widehat M,\widehat h)$) be the oriented Riemannian manifold obtained from $(M, g)$ (resp. $(M, h)$) by metric surgery along $\cS$. By construction, each of the metrics $\widehat g, \widehat h$ has scalar curvature greater than or equal to $1$ and geometry bounded by some constant $C$. 

Since $\cS$ is a spherical splitting, every component of $\widehat M$ is closed. We denote by $\{M_\ell\}$ the collection of those components. Applying Theorem~\ref{thm:marques deforming2} to each $M_\ell$, we get a family of continuous paths $\widehat g_t^{(\ell)} $ in $\cR_1(\widehat M_\ell)$ and positive diffeomorphisms $\widehat \psi^{(\ell)} : \widehat M_\ell \to \widehat M_\ell$ such that for each $l$, we have $\widehat g_0^{(\ell)}=\widehat g_{\vert M_\ell}$ and $\widehat \psi^{(\ell)}_{\ast}\widehat g_1^{(\ell)} = \widehat h_{\vert M_\ell}$.

Applying Lemma~\ref{lem:combining paths} and Remark~\ref{rem:combining paths}, we  get a path $\widehat h_t$ in $\cR_1^{bg}(\widehat M)$. Denote by $p_\alpha^\pm$ the tips in $\widehat M$ of the caps $\cS^\pm_\alpha$ added by the metric surgery process. Note that the tips are the same with respect to both metrics $\widehat g$ and $ \widehat h$, as well as the balls $B_{2\rho_C}(p_\alpha^\pm)$ for a suitable radius $\rho_C$; furthermore, these balls are round and pairwise disjoint.

Denote by $\widehat g_\#=(\widehat g)_\#(\{p_\alpha^\pm\})$ (resp.~$\widehat h_\#=(\widehat h)_\#(\{p_\alpha^\pm\})$) the GL-sum metric constructed with GL-parameters $(\rho_C,1,\eta_C)$ for an appropriate $\eta_C$. We identify the corresponding manifold $\widehat M_\#$ with $M$ as in the last paragraph of Section \ref{sec:cn metric surgery}. 

We deduce from Lemma \ref{lem:gl-surg} that  $\widehat g_\#$ is isotopic to $g$ and that $\widehat h_\#$ is isotopic to $h$ (with uniformly bounded geometry). To finish the proof, we need to isotope $\widehat g_\#$ to $\widehat h_\#$ modulo diffeomorphism.

Let $\widehat \psi : \widehat M \to \widehat M$ be the positive diffeomorphism defined by setting $\widehat \psi = 	\widehat \psi^{(\ell)}$ on $\widehat M_\ell$. Define  $q_\alpha^\pm = {\widehat \psi}^{-1}(p_\alpha^\pm)$. For each $\alpha$ and each  $\epsilon\in\{1,-1\}$, choose a positive orthonormal basis  $\{e_{k,\alpha}^\epsilon\}$ at $p_\alpha^\epsilon$.  Set $f_{k,\alpha}^{\pm}=\widehat{\psi}^\ast(e_{k,\alpha}^\pm)$. 
	Since $\widehat \psi$ is a positive isometry from $(\widehat M,\widehat{g_1})$ to $(\widehat M, \widehat h)$, sending $\{q_\alpha^\pm\}$ to $\{p_\alpha^\pm\}$ and $\{f_{k,\alpha}^{\pm}\}$ to $\{e_{k,\alpha}^\pm\}$, it induces an isometry from $(\widehat M,\widehat{g_1})_\#(\{q_\alpha^\pm\}),\{f_{k,\alpha}^{\pm}\})$ to $(\widehat M, \widehat h)_\#(\{p_\alpha^\pm\},\{e_{k,\alpha}^{\pm}\}))=(\widehat M_\#,\widehat h_\#)$. Thanks to Lemma~\ref{lem:som gl gb}, we can isotope the former manifold---modulo diffeomorphism---to  $(\widehat M,\widehat g_0)_\#(\{p_\alpha^\pm\})=(M_\#,g_\#)$. This completes the proof of Lemma~\ref{lem:same cS isotopic}.
\end{proof}

\subsection{Proof of Proposition~\ref{prop:isotopy-straight2}}
\label{subsec:isotopy-straight2}

We start with a technical lemma.
\begin{lem}\label{lem:pasting metrics} Let $M$ be a $3$-manifold. Let $(g,\cS, g', \cS')$ be a quadruple such that 
	\begin{enumerate}[\hspace{0.2cm} (i)]
		\item  $g \in \cR_1^{bg}(M)$ is a GL-metric and 
		$\cS=\{S_\alpha \}$ is a spherical splitting straight for $g$,
		\item  $g' \in \cR_1^{bg}(M)$ and $\cS'=\{S'_\beta\}$ is a spherical system straight for $g'$, 
		\item Components of $\cS$ and $\cS'$ are mutually  disjoint.
	\end{enumerate}
	Then there exists $h=h(g,\cS, g',\cS') \in \cR_1^{bg}(M)$ such that  $h=g$ near $\cS$ and $h=g'$ near $\cS'$.
\end{lem}

\begin{proof} 
By Theorem~\ref{thm:BBM11}, there exists a finite collection $\cF$ of spherical manifolds such that $M$ is a connected sum of members of $\cF$. 

For each $\alpha,\beta$,	let $U_\alpha$ (resp. $U'_\beta $) be a straight tube w.r.t.~$g$ containing $S_\alpha $ (resp.~w.r.t.~$g'$ containing $S'_\beta$). Assume that all theses	tubes are pairwise disjoint. Let $\widehat M$ be a $3$-manifold obtained from $M$, by splitting every $U_\alpha$ and $U'_\beta$ along a pair of straight spheres on each side of the $S_\alpha$ and of the $S'_\beta$, and glueing $3$-balls to the boundary. All components of $\widehat M$ are closed. All components containing some $S_\alpha$ or some $S'_\beta$ are topological $3$-spheres. All other components are connected sums of members of $\cF$.

We define a metric $\widehat h$ on $\widehat M$ as follows. On components containing 
	some $S_\alpha$ (resp. $S'_\beta$) we let $\widehat h$ be the metric obtained from $g$ (resp. $g'$) by metric surgery on the straight tubes. In particular  $\widehat h=g$ near $S_\alpha$ and $\widehat h=g'$ near $S'_\beta$. On the other components we let $\widehat h$ be a GL-sum of round metrics on members  of $\cF$, of uniformly bounded geometry and scalar curvature $\geq 2$. 
	
Finally we define $h\in \cR_1^{bg}(M)$ by doing the GL-sum of $(\widehat M,\widehat h)$, in such a way that $h=g$ near $S_\alpha$ and $h=g'$ near $S'_\beta$.
\end{proof}

\begin{proof}[Proof of Proposition \ref{prop:isotopy-straight2}]
Let $\cS=\{S_\alpha\}$ be a spherical splitting that is straight for $g$. By local finiteness of $\cS$ and $\cS'$, up to replacing $\cS$ by a sub-system there exists a sub-system $\cS''\subset \cS'$ such that every component of $\cS$ is disjoint from every component of $\cS''$.
	
	By Lemma~\ref{lem:pasting metrics}  there exists a metric 
	$h'=h(g,\cS,g', \cS' ) \in \cR_1^{bg}(M)$ such that $h'=g$ near $\cS$ and $h'=g'$ near $\cS''$. Applying Lemma~\ref{lem:same cS isotopic} to the triple $(h',g,\cS)$ we see that $h'$ and $g$ are isotopic modulo diffeomorphism.

	 By Lemma~\ref{lem:pasting metrics}  there exists a metric $h=h(h',\cS'',g', \cS' \setminus \cS'') \in \cR_1^{bg}(M)$ such that  $h=h'$ near $\cS'' $ and $h=g'$ near $\cS'\setminus \cS''$. From the fact that $h'=g'$ near $\cS''$, we deduce that $h=g'$ near $\cS'$. It follows from Lemma~\ref{lem:same cS isotopic} applied to the triple $(h,h',\cS'')$ that $h$ is isotopic to $h'$ modulo diffeomorphism. Hence $h$ is isotopic to $g$ modulo diffeomorphism, which is the required conclusion.
\end{proof}

\section{Surgical solutions of Ricci flow} \label{sec:surgical_solution}

In this section we recall the basic properties of the surgical solutions constructed in \cite{B2M:scalar}. 

\subsection{Evolving metrics and surgical solutions}

\begin{defi}
An \emph{evolving Riemannian 3-manifold} is a family $\{(M(t),g(t))\}_{t \in I}$ of (possibly empty or disconnected) Riemannian $3$-manifolds indexed by an interval $I\subset\RR$. It is \emph{piecewise $C^1$-smooth} if there exists a subset $J$ of $I$ which is discrete as a subset of $\RR$ and satisfies the following conditions:
\begin{enumerate}[\indent (i)]
 \item On each connected component of $I \setminus J$, $t \mapsto M(t)$ is constant and $t\mapsto g(t)$ is $C^1$-smooth.
 \item For each $t_0 \in J$, $M(t)=M(t_0)$ for all $t <t_0$ close enough to $t_0$ and $t \mapsto g(t)$ is left continuous at $t_0$.
 \item For each $t_0 \in J \setminus \sup{I}$, $t \mapsto (M(t),g(t))$ has a right limit at $t_0$, denoted by $(M_+(t_0),g_+(t_0))$.
\end{enumerate}
A time $t \in I$ is \emph{singular} if $t\in J$ and regular otherwise.
\end{defi}

\begin{defi} A piecewise $C^1$-smooth evolving $3$-manifold $(M(t),g(t))$ is called a \emph{surgical solution} if the following holds:
\begin{enumerate}[\indent (i)]
 \item The Ricci Flow equation $\frac{dg}{dt}=-2\Ric_{g(t)} $ is satisfied at all regular times.
 \item For each singular time $t$ we have $\Rmin(g_+(t)) \geq \Rmin(g(t))$.
 \item For each singular time $t$ there is a locally finite collection $\cS(t)$ of disjoint embedded $2$-spheres in $M(t)$ and a 
 manifold $M'(t)$ such that 
 \begin{enumerate}[\indent (a)]
 \item $M'(t)$ is obtained from $M(t) \setminus \cS(t)$ by capping-off $3$-balls;
 \item $M_+(t)$ is a union of connected components of $M'(t)$ and $g_+(t)=g(t)$ on $M_+(t) \cap M(t)$;
 \item Each component of $M'(t) \setminus M_+(t)$ is spherical or diffeomorphic to $\RR^3$, $S^2\times S^1$, $S^2 \times \RR$, $\RR P^3 \#\RR P^3$ or a punctured 
 $\RR P^3$.
 \end{enumerate}
\end{enumerate}
For a singular time $t$, the components of $M'(t) \setminus M_+(t)$  are called the \bydef{discarded components}. The surgical solution is \bydef{extinct} if for some $t$, all components are discarded, i.e.~$M_+(t)=M(t)=\emptyset$.
\end{defi}

\begin{rem}\label{rem:extinct}
When the initial metric has uniformly positive scalar curvature and each $g(t)$ is complete with bounded sectional curvature, we deduce from the maximum principle and Property~(ii) a lower bound for the scalar curvature which blows up in finite time; thus under these hypotheses the solution must be extinct.
\end{rem}

 In \cite{B2M:scalar} we constructed a special kind of surgical solutions, called $(r,\delta,\kappa)$-surgical solutions, which use the same three parameters as Perelman's Ricci flow with surgery and one more which serves as a curvature threshold to trigger the surgery. We will not need the precise definition here. For us, the two main properties of $(r,\delta,\kappa)$-solutions are that whenever $t$ is a singular time, the manifold $(M_+(t),g_+(t))$ is obtained from $(M(t),g(t))$ by metric surgery on $\epsi$-necks, and the discarded components are covered by so-called canonical neighbourhoods. We have already reviewed metric surgery in Section~\ref{sec:cn metric surgery}. In the next subsection, we discuss the notion of a canonical neighbourhood.

\subsection{Canonical neighbourhoods and locally canonical metrics}\label{subsec:cn}
 Let $(M,g)$ be a Riemannian 3-manifold and $\epsi,C>0$ be constants. There are four types of canonical neighbourhoods: necks, caps, $\epsi$-round components and $C$-components. We already defined the notion of an $\epsi$-neck in Section~\ref{sec:cn metric surgery}. A component $X$ of $M$ is \emph{$\epsi$-round} if after scaling to make $R(x)=1$ at some point, $X$ is $\epsi$-close to a round metric of scalar curvature one. A component of $M$ is a \emph{$C$-component} if it is diffeomorphic to $S^3$ or $\textbf{R}P^3$, and has positive sectional  curvature and geometry bounded by $C$ after scaling. More precision can be given on the geometry of these neighbourhoods 
 (see~e.g.~Definition 4.2.8 in \cite{B3MP} or  \cite{marques:deforming} p.~837). 
 
 An \emph{$\epsi$-cap} is an open subset $\cC \subset M$ diffeomorphic to a $3$-ball or to $\textbf{R}P^3$ minus a $3$-ball, with an $\epsi$-neck $N \subset \cC$ such that 
 $\overline{Y}=\cC - N$ is a compact submanifold with boundary. The boundary $\partial \overline{Y}$ of the core $Y$ (interior of $\cC - N$) is required to be the middle sphere of an $\epsi$-neck. An $(\epsi,C)$-cap is an $\epsi$-cap such that, after rescaling so that $R(x)=1$ for some point $x$ in the cap, the diameter, volume and curvature ratios at any two points are bounded by $C$.

\begin{defi}
A point $x$ in $(M,g)$ is said to be centre of an \bydef{$(\epsi,C)$-canonical neighbourhood} if it is centre of an $\epsi$-neck, or centre of an $(\epsi,C)$-cap, or is contained in an $\epsi$-round component or a $C$-component.  If each point of $(M,g)$ is centre of a $(\epsi,C)$-canonical neighbourhood, we will say that $g$ is  \emph{$\epsi$-locally canonical}. 
\end{defi}

We now fix the constants $\epsi,C$, refining the choice made in \cite{B2M:scalar}, so that the interpolation lemmas of~\cite{marques:deforming} hold.  Set
 $\epsi = \min(\delta_0,\epsi_3)$, where $\delta_0$ is the constant from
 Theorem~\ref{thm:surgery} and $\epsi_3$ is from \cite[Lemma 6.3]{marques:deforming}. Then set $C=\max(100,2C_{\mathrm{sol}}(\epsi/2),2C_{\mathrm{st}}(\beta(\epsi)\epsi/2))$ as in~\cite[p.~947]{B2M:scalar}.

From \cite[section 5.2]{B2M:scalar} we gather the following existence result, which can be taken as a black box for the rest of the paper: 
\begin{theo} \label{thm:existence_surgical} Given a positive number $A$, there exist a positive number $\tau$ and a tuple $\mathcal{Q}=\mathcal{Q}_{[\epsi^{-1}]}$ with the following property. 
Let $(M,g)$ be a complete Riemannian manifold with geometry bounded by $A$, and with scalar curvature greater than or equal to $1$. 
Then there exists an extinct surgical solution $(M(\cdot),g(\cdot))$ defined 
on $[0,2]$ such that $(M(0),g(0))= (M,g)$ and satisfying the following properties: 
\begin{enumerate}
 \item For all $t\in [0, 2]$, we have $g(t) \in \mathcal{R}^{Q_0}_1(M(t))$,
 \item The solution is smooth on $[0, \tau ]$ and, for all $t\in [\tau, 2]$, we have $g(t) \in \mathcal{R}^{\mathcal{Q}}_1(M(t))$.
 \item At each singular time $t$, the manifold $M'(t)$ has a metric $g'(t)\in \mathcal{R}_1^{\mathcal Q}(M'(t))$ obtained from $g(t)$ by metric surgeries in $\epsi$-necks, and such that
   \begin{enumerate}
     \item $g_+(t)=g'(t)$ on $M_+(t)$.
     \item All discarded components of $(M'(t),g'(t))$ are $\epsi$-locally canonical.
  \end{enumerate}
 \item The number of singular times is finite and bounded by a constant depending on $A$.
\end{enumerate}
\end{theo}

In the proof of Theorem~\ref{thm:isotopy_to_GL}, Theorem~\ref{thm:existence_surgical} will allow to isotope an arbitrary metric to a GL-sum of $\epsi$-locally canonical ones. Thus it is useful to be able to isotope locally canonical metrics to GL-metrics. This is the purpose of the next lemma:
\begin{lem}\label{lem:loccan_to_GL}
Given $\cA_{[\epsi^{-1}]}$ there exists $\cB_{[\epsi^{-1}]}$ such that the following holds. Let $M$ be a connected $3$-manifold and let $g$ be an $\epsi$-locally canonical metric belonging to  $\cR_1^{\cA_{[\epsi^{-1}]}}(M)$. Then there exists an isotopy $ g_t \in \cR_1^{\cB_{[\epsi^{-1}]}}(M)$ such that $g_0=g$  and $g_1$ is a GL-metric.   
\end{lem}

\begin{proof}
If $M$ is compact, every metric is a GL-metric (the empty collection is a spherical splitting), so we assume that $M$ is non-compact. Let $g \in \cR_1^{\cA_{[\epsi^{-1}]}}(M)$ be an $\epsi$-locally canonical metric. From~\cite[proof of Proposition~7.2]{B2M:scalar} we have three cases:
\begin{enumerate}
     \item $M$ is covered by $\epsi$-necks and is diffeomorphic to $S^2\times \RR$, or
     \item $M$ is covered by $\epsi$-necks  and one $\epsi$-cap diffeomorphic to $B^3$, and $M$ is diffeomorphic to $\RR^3$, or
     \item $M$ is covered by $\epsi$-necks  and one $\epsi$-cap diffeomorphic to $\RR P^3\setminus \{\textrm{point}\}$, and $M$ is diffeomorphic to $\RR P^3\setminus \{\textrm{point}\}$.
\end{enumerate}

Take a maximal family of disjoint $\epsi$-necks. From the above we see that the middle spheres of these necks form a spherical splitting of $M$. Thus it suffices to deform the metric $g$ to a metric with respect to which these spheres are straight.

In each neck we apply the following deformation. We use the parametrisation to identify a neighbourhood of the middle sphere with $S^2\times (-4, 4)$. Let $h$ be the pulled-back rescaled metric on $S^2\times (-4, 4)$. Now, we choose a cutoff function $\eta : (-4,4) \to [0,1]$ such that $\eta(s)=1$ on $[-2,2]$ and $\eta(s)=0$ when $|s|$ is close to $4$. Then for $t\in [0,1]$ we set $h_t = (1-\eta(s)t)h+\eta(s)t\gcyl$, where $s$ is the radial coordinate. We observe that $h_1= \gcyl$ on $S^2\times [-2,2]$ and that the deformation $h_t$ is constant near the boundary of $S^2\times (-4, 4)$. It follows that the rescaled pushed-forward metric on $M$ defines the required isotopy.
 \end{proof}

\section{Isotopies of uniformly bounded geometry} \label{sec:isotopy bounds}

In this section we show that if $(M,g)$ is isotopic to a GL-metric in some $\cR_\sigma^{\cA}$ then, after a GL-sum for suitable parameters,  $(M_\#,g_\#)$ is isotopic to a GL-metric in some 
$\cR_{\sigma'}^{\cB}$, with $\cB$ and $\sigma'$ depending on $\cA$, $\sigma$ and the parameters. The following proposition gives a more precise statement.

\begin{prop} \label{prop:sum_isotopies} Let $\sigma>0$.  For all $\cA=\cA_{[\epsi^{-1}]}$ and $\rho \in (0, \min\{c(3){A_0}^{-1/2},1\})$, there exists $\bar{\cA}^{\#}=\bar{\cA}^{\#}_{[\epsi^{-1}]}$ and $\eta>0$ with the following property. Let $(M,g)$ be an oriented $3$-manifold and $\{ (p_\alpha^-, p_\alpha^+)\}$ be a family of pairs of points. Set $\cP=\{p_\alpha^\pm\} \subset M$ and $(M_\#,g_\#)=(M,g)_\#(\{(\rho,\sigma,\eta)\},\{p_\alpha^\pm\})$. Suppose that the following assertions hold:
	\begin{enumerate}[\hspace{0.3cm} (a)]
		\item For all $x, y \in \cP$, if $x \not= y $, then $B_{g}(x,3\rho) \cap B_g(y,3\rho)=\emptyset$, and $g$ is round on $B_{g}(x,3\rho)$.
		\item The metric $g$ is isotopic in $\cR_\sigma^{\cA}(M)$ to some GL-metric.
	\end{enumerate}
	 Then $g_\#$ is isotopic in $\cR_{9\sigma/10}^{\bar{\cA}^{\#}}(M_\#)$ to some GL-metric.
\end{prop}

We first explain some of the key ideas informally. Let $g_t$ be an isotopy in $\cR_\sigma^{\cA}(M)$ from $g$ to some GL-metric. We would like to apply Proposition~\ref{prop:continuousGL},  moving base points $p_\alpha^\pm$ along continuous paths if necessary, to get a continuous GL-sum $(M,g_t)_\#$.  The main difficulty is that we have no control on $d_{g_y}(x,y)$, for $x\not= y\in \cP$. For example, a compact component of $M$ containing  an unspecified number of $p_\alpha^\pm$ could have a diameter becoming small along the path of metrics $g_t$; therefore the points  $p_\alpha^\pm$  become close to each other. Moreover, a reparametrisation trick, as in Lemma \ref{lem:som gl gb}, is not possible as we deal with possibly non-compact components and the metric $g_1$ does not a priori verify the assumption a).

The trick is to first deform continuously $(M_\#,g_\#)$ into a manifold in which is obtained by connecting sum around points which are far away. Looking at $p_\alpha^+$ and $p_\alpha^-$, the GL-sum gives rise to a metric on the cylinder  joining neighbourhoods of these two points which is close to $\gcyl$. We modify it continuously, applying twice Lemma \ref{lem:conf-cyl}, so that it becomes $ B_\rho(p_\alpha^-)\, \# \, S^3(p_\alpha^-)\, \# \, S^3(p_\alpha^+) \, \# \, B_\rho(p^+)$ where $S^3(p_\alpha^\pm)$ are round $3$-spheres which will undergo only very small deformations during the process (see Sublemma \ref{sublemma:inflate}). Then we perform a surgery on this topological cylinder by cutting it between the two spheres and glueing $3$-balls on either side.

Doing this for every $\alpha$ creates some distance between the points $p_\alpha^+$ and $p_\alpha^-$, so that we can deal with separately. Note that they may remain in the same connected component or not. The next step is to produce a set $\cP_{sep} \subset \cP$ such that its points remain sufficiently separated during the isotopy (see Lemma \ref{lem:divide P}). Each $p_0\in \cP_{sep}$ has a neighbourhood which contains only finitely many points of $\cP$, say $(p_1, \dots, p_n)$. We then do the GL-sum of $S^3(p_0)$ with $n$ round spheres and we continuously move each $p_i$ in one of them. They will be only slightly deformed and hence the points will remain far apart.

This process deforms $(M_\#,g_\#)$ into a manifold isometric to a GL-sum $(M \sqcup X,g \sqcup g_X)_\#$, where $g_X$ is a GL-sum of round spheres, and where the connected sum is made between $M$ and $X$ along a sparse set of points, which remain far from each other along $g_t \sqcup g_X$. We can then safely apply Proposition~\ref{prop:continuousGL} to get a continuous path $(g_t \sqcup g_X)_\#$ with controlled geometry. We summarise the idea saying that it amounts to 'externalising' the connected sum. 

We now turn to the formal proof.

\begin{proof}[Proof of Proposition~\ref{prop:sum_isotopies}]

	Let $\cA$, $\rho$ and $\sigma$ be as in the statement of Proposition \ref{prop:sum_isotopies} and let $g_t$ be an isotopy in $\cR_\sigma^{\cA}(M)$ such that $g_0=g$ and $g_1$ is a GL-metric.
	
	Let  $\eta=\eta(\cA_2)$ be the constant given by Lemma \ref{lem:eta}. Let $\cA'=\cA^{\#}(\cA,\rho,\sigma,\eta)$ be the constant of Proposition \ref{prop:continuousGL}. Set $\sigma'=\frac{9}{10}\sigma$.  Then  $(M_\#,g_\#)=(M,g)_\#(\{(\rho,\sigma,\eta)\},\{p_\alpha^\pm\})$  has scalar curvature greater than $\sigma'$ and geometry bounded by $\cA'$. Let $\rho'=\min\{c(3){\cA'_0}^{-1/2},\frac{1}{10}\rho\}$, $\eta'=\eta(\cA'_2)$ be the constant given by Lemma \ref{lem:eta}, $\bar \sigma$ and $\bar \cA$ be the constants given by Lemma \ref{lem:conf-gamma} applied to $\gamma_1=\gamma(\rho,\sigma,\eta)$, $\gamma_2=\gamma(\rho',\sigma',\eta')$. We denote by $\bar{\bar \sigma}=\bar{\bar \sigma}(\rho,\sigma,\eta)$ and $\bar{\bar{\cA}}=\bar{\bar{\cA}}(\rho,\sigma,\eta)$ the constants given by Lemma~\ref{lem:conf-cyl}.
	We then define $\cA'' = \min\{\cA^{\#}(\cA', \rho',\sigma'), \bar \cA, \bar{\bar{\cA}}\}$ and $\sigma''=\min\{(9/10)^2\sigma,\bar{\bar \sigma}\}$.

	\begin{defi} \label{def:separated}
		We say that 
		a subset $\cP \subset M$ is \emph{$\rho$-separated for an isotopy $g_t$} if 
		$$ \forall x \not= y \in \cP, \quad \bigcup_t B_{g_t}(x,\rho) \cap  \bigcup_t B_{g_t}(y,\rho) = \emptyset.$$ 
	\end{defi}
	In particular, if $\cP=\{p_\alpha^\pm\}$ is $3\rho$-separated, then there exist pairwise disjoint topological $3$-balls $W_\alpha^\pm$ such that $B_{g_t}(p_\alpha^\pm,2\rho) \subset W_\alpha^\pm$ for all $t$, as in the assumptions of Proposition \ref{prop:continuousGL} (with fixed $p_\alpha^\pm$).

	The first step is to divide points of $\cP$ into disjoint open connected subsets $U_j$ of $M$, such that each $U_j$ contains finitely many elements of $\cP$ and such that they are in some sense not to close one from each other. This will allow,  by moving the base points independently in each $U_i$, to deform the GL-sum $g_\#$ with controlled geometry.   
				
			\begin{lem}\label{lem:divide P} Up to moving points of $\cP$ by a distance less than $20\rho' \leq 2\rho$, there exist a subset $\cP_{sep} \subset \cP$, a family $(U_j)$, $j\in J\subset \textbf{N}$,  of disjoint open connected subsets of $M$ and a family $(A_k)$ of open subsets of $M$, such that,
				\begin{enumerate}[\hspace{0.3cm} $\bullet$]
					\item $\cP_{sep}$ is $1$-separated for $g_t$, for all $t$.
					\item The $7\rho'$-neighborhood of $\cP$ is contained in $\bigcup\limits_{j\in J} U_{j}$ 
 and for all $j$, $|\cP \cap U_j| < \infty$ and $|\cP_{sep} \cap U_j| = 1$.
					\item $\bigcup\limits_{k} A_k = \bigcup\limits_{j} U_j$ and each $A_k$ is a union of finitely many $U_j$.
					\item For all $k\not=l$, $d_g(A_k,A_l) \geq 6\rho'$.\\
				\end{enumerate}
			\end{lem}

\begin{proof}		
Without loss of generality, we assume that $M$ is connected. If $\cP$ is finite, then we set $U_0=M$, $A_0=U_0$ and $\cP_{sep}=\{p_0\}$, where $p_0$ is an arbitrary point of $\cP$.

Suppose that $\cP$ is infinite. Up to moving slightly some points of $\cP$, we will define $\cP_{sep}=\{p_0,p_1,\ldots,\} \subset \cP$ and $(U_j)_{j \in \NN}$ with $\cP_{sep} \cap U_j=p_j$.  We use the following property, which is easily derived from the continuity of $t \mapsto g_t$ and the compactness 
			of $[0,1]$: {\it for every $p \in M$ and $R_0>0$, there exist a compact subset $K \subset M$ and $R_1>0$ such that 
				$$ \forall s,t \in [0,1],\quad  B_{g_s}(p,R_0) \subset K  \subset B_{g_t}(p,R_1).$$}
			
			This property implies that for every $p\in M$, for every sufficiently large $R$, for every $x \in M$, if $d_g(p,x) \geq R$, then $\{p,x\}$ is $1$-separated for $g_t$. Let us fix some $p=p_0 \in \cP$ and a radius 
			$R_0>1$. Set $\tilde R_0=R_0 -20\rho'$. Up to moving points of $(\overline B_g(p_0,R_0) \setminus B_g(p_0,\tilde R_0)) \cap \cP$ inside $B_g(p_0,\tilde R_0)$, we can assume that the closed annulus $\overline B_g(p_0,R_0) \setminus B_g(p_0,\tilde R_0)$ does not intersect $\cP$. Let $\tilde U_0= B_g(p_0,\tilde R_0+7\rho')$.  We add to $\tilde U_0$ all the connected components of $M \setminus \tilde U_0$ containing only finitely many elements of $\cP$; there are finitely many such components.  This defines an open connected subset $U_0 \subset M$ such that  
			\begin{itemize}
				\item The set $U_0 \cap \cP$ is finite.
				\item For every component $C$ of $M \setminus U_0$, the set $C \cap \cP$ is infinite.
				\item We have $\partial U_0 \subset \partial B(p_0,\tilde R_0+7\rho')$.
				\item The $7\rho'$-neighborhood of $\cP \cap U_0$ is contained in $U_0$. 
			\end{itemize}
			
			Set $A_0=U_0$. 
			Let $C_1,\ldots,C_k$ be the connected components of $M \setminus (B(p_0,R_0-7\rho') \cup A_0)$. Note that $d(C_i, A_0)\geq 6\rho '$ and that the $7\rho'$-neighborhood of $\partial C_i$ does not intersect $\cP$. Choose for each $i\in \{1,\ldots,k\}$ a point $p_i \in C_i$ far enough from $\partial C_i$ and such that for every $t\in [0,1]$, the set $\{p_0,p_1,\ldots,p_k\}$ is $1$-separated w.r.t.~$g_t$. Choose  a radius $R_1$ which is bigger than $2d(p_i,\partial C_i)$ for all $i\ge 1$. Set $\tilde R_1=R_1 - 20\rho'$, and for every $i$, set $\tilde U_i:=\Int(C_i) \cap B(p_i,\tilde R_1+7\rho')$.	Perform, in restriction to $\Int(C_i)$, the previous construction starting from $\Int(C_i) \cap B(p_i,R_1)$, i.e.~move points of $\Int(C_i) \cap B(p_i,R_1) \cap \cP$ inside $\Int(C_i) \cap B(p_i,\tilde R_1)$.
 
Define $U_i \subset \Int(C_i)$ by adding to $\tilde U_i$ connected components of $\Int(C_i) \setminus \tilde U_i$ containing finitely many elements of $\cP$. Thus $U_i \cap \cP$ is finite, and for each component $C$ of $\Int(C_i) \setminus U_i$, the set $C \cap \cP$ is infinite. Furthermore, $\partial U_i \cap \Int (C_i)$ is contained in $\partial B(p_i,\tilde R_1+7\rho')$, and the $7\rho'$-negihborhood of 
$\cP \cap U_i$ is contained in $U_i$.

Observe that the $U_i$'s are pairwise disjoint. Setting $A_1=U_1\cup\dots\cup U_k$, we have $d(A_0,A_1)\geq 6\rho '$. We now iterate the construction on $M\setminus (A_0\cup A_1\cup B(p_0, R_0-7\rho')\cup \dots \cup B(p_i, R_1-7\rho'))$.
				\end{proof}
	
	Up to deforming slightly $g_\#$ by moving the base points, we can assume that $\cP \supset \cP_{set}$ as given by Lemma \ref{lem:divide P}. The deformation of $g_\#$ which `externalises' most of its GL-sums is given by the following lemma.

	\begin{lem}\label{lem:deform gl} There exist a $3$-manifold $X$, a GL-metric $g_X \in \cR_1^{\cA}(X)$, and an injective map $x : \cP_{sep} \to X$ with the following properties:
		\begin{enumerate}
     		\item $x(\cP_{sep}) \subset X$ is $3\rho$-separated for $g_X$, and $g_X$ is round on its $3\rho$-neighbourhood.  
			\item The GL-sum $(M \sqcup X, g \sqcup g_X)_\#(\{(\rho,\sigma,\eta)\}, \{(p,x(p)) \mid p \in \cP_{sep}\})$ has scalar curvature greater than $\sigma'$ and geometry bounded by $\cA'$, and $(M \sqcup X)_\#$ is diffeomorphic to $M_\#$. 
			\item $(M_\#,g_\#)$ is isotopic modulo diffeomorphism to $(M \sqcup X, g \sqcup g_X)_\#$ through metrics with scalar curvature greater than $\sigma''$ and geometry bounded by $\cA''$.
		\end{enumerate}
	\end{lem}

\begin{proof}	
We construct an isotopy $\{g_\#^t\}_{t \in [0,5]}$ on $M_\#$, and a manifold $(X,g_X)$ as above, such that $(M_\#,g_\#^5)$ is isometric to  $(M \sqcup X, g \sqcup g_X)_\#$. 

Recall that the GL-sum construction associates to each pair $p^\pm \in \{p_\alpha^\pm\}$ a GL-neck $$B_\rho(p^-)\# B_\rho(p^+) \approx S^2 \times I \subset M_\#$$ where the metric $g_\#$ is cylindrical near the middle of the neck. For each $p^\pm$ we denote by $S^3(p^-)$ and $S^3(p^+)$ two round $3$-spheres with scalar curvature $\sigma$, and we fix two points $x(p^-) \in S^3(p^-)$ and $x(p^+) \in S^3(p^+)$.

	\begin{slem}\label{sublemma:inflate} There exists an isotopy $\{g_\#^t\}_{t \in [0,1]}$, which deform each GL-neck $B_\rho(p^-)\# B_\rho(p^+)$, without changing the metric near the boundary, through metrics with scalar curvature $\geq \bar{\bar \sigma}$ and geometry bounded by $\bar{\bar \cA}$,  into the GL-sum 
			$$ B_\rho(p^-)\, \# \, S^3(p^-)\, \# \, S^3(p^+) \, \# \, B_\rho(p^+)$$
			made  with parameters $(\rho,\sigma,\eta)$ at  $\{ (p^-,x(p^-)), (-x(p^-),-x(p^+)), (x(p^+),p^+)\}$.
	\end{slem}		
		
	\begin{proof} Apply Lemma \ref{lem:conf-cyl} twice.
	\end{proof}
	
Undoing each GL-sum $S^3(p^-) \# S^3(p^+)$ of $(M_\#, g_\#^1)$ splits the manifold into a manifold $(\widehat M, \widehat g_1)$ where $\widehat M$ is diffeomorphic to $M$ and $\widehat g_1 \in \cR_{\sigma'}^{\cA'}(\widehat M)$. Our goal is to construct an isotopy $\{\widehat g_t\}_{t \in [1,5]}$ on $\widehat M$. The isotopy $g_\#^t$ will be obtained from $\widehat g_t$ by reconnecting each pair $S^3(p^-) \# S^3(p^+)$.
	
	Observe that $(\widehat M,\widehat g_1)$ is simply the manifold obtained from $(M,g)$ by replacing each ball $B_\rho(p)$ by $B_\rho(p)\, \# \, S^3(p)$, for all $p \in \{p_\alpha^\pm\}$. We denote by $H(p) \subset S^3(p)$ the hemisphere centred on $-x(p)$ (which is on the opposite of the hemisphere connected to $B_\rho(p)$). Therefore $H(p) \subset B_\rho(p)\, \# \, S^3(p) \subset \widehat M$. 
			
		\begin{slem} \label{sublemma:deux} There exist a $3$-manifold $\widehat X$ and a GL-metric $\widehat g_X \in \cR_{\sigma'}^{\cA'}(\widehat X)$, an injective map $\widehat x:\cP_{sep} \to \widehat X$ and an isotopy $(\widehat g_t)_{t \in [1,5]}$ on $\widehat M$ such that:
			\begin{enumerate}
				\item $\cP_{sep}$ is $3\rho$-separated for  $(\widehat g_t)$, and  $\widehat x(\cP_{sep})$ is $3\rho$-separated for $g_{\widehat X}$.
				\item $\widehat g_t = \widehat g_1$ on all hemispheres $H(p^-)$, and $H(p^+)$.
				\item $\widehat g_t$, resp.  $\widehat g_5$, has scalar curvature greater than $\sigma''$ (resp.~$\sigma'$) and geometry bounded by $\cA''$ (resp.~$\cA'$.)
				\item $(\widehat M,\widehat g_5)$ is isometric to $(M \sqcup \widehat X, g \sqcup g_{\widehat X})_\#(\{(p,\widehat x(p)) \mid p \in \cP_{sep}\})$. Moreover $H(p) \subset \widehat X$ for all $p \in \cP$. 
			\end{enumerate}
		\end{slem}

 	\begin{proof} 
 		 By Lemma \ref{lem:divide P}, the set $\cP$ is divided into open connected subsets $U_j$ of $M$ and these subsets are gathered into open subsets $A_k$ such that each $A_k$ is a union of finitely many $U_i$ and $d(A_k, A_l)\geq 6\rho '$ when $k\not = l$. Let $V_k \subset (M,g)$ denote the $3\rho'$-neighbourhood of $A_k$. By construction, the $V_k$'s are pairwise disjoint. For each $k$, let  $\widehat V_k$ be the corresponding subset in $\widehat M$. The $\widehat V_k$'s are also pairwise disjoint.

 		  We will contruct isotopies $(\widehat g_t^{\widehat V_k})_{t \in [1,5]}$ on $\widehat M$, with support in $\widehat V_k$, and then define the isotopy $(\widehat g_t)_{t \in [1,5]}$ on $\widehat M$ by as $$ \left\{\begin{array}{cccl}
 		 		\widehat g_t & = & \widehat g_t^{\widehat V_k} & \textrm{ on } \widehat V_k,\, \forall k\\ 
 		 		\widehat g_t & = & \widehat g_1 & \textrm{ on } \widehat M \setminus \bigcup_k \widehat V_k.
 		 		\end{array}  \right.$$

 		 The $(\widehat g_t^{\widehat V_k})$ are constructed independently of each other. Let us fix $k$ and the subset $\widehat V_k$. For the sake of simplicity we now drop the index $k$. To $\widehat V= \widehat V_k$ we associate a compact GL-manifold $(\widehat X_{\widehat V},g_{\widehat X_{\widehat V}})$, which in fact is simply a (possibly non connected) GL-sum of a finite number of round $3$-spheres, and an isotopy $\widehat g_t^{\widehat V}$ defined on $\widehat M$ with support in $\widehat V$, starting from $\widehat g_1$. The set $V=V_k$ is the $3\rho '$-neighbourhood of  $A=A_k$. The open set $A$ is a union  $U_{i_1}\cup\dots\cup U_{i_N}$.

For simplicity, we first discuss the case where $N=1$ and perform a four-step deformation of $\widehat g_1$ into a metric $\widehat g_5$. The general case will be tackled later.

 Droping indices again, we consider $U=U_{i_1}$ and denote by $\cP\cap U=\{ p_1, \dots , p_n\}$ where $\cP_{sep}=\{ p_1 \}$.

		Step 1.{\it Deform $(B_\rho(p_1) \# S^3(p_1), \widehat g_1)$, without changing the metric near the boundary nor on the hemisphere $H(p_1)$, through metrics with scalar curvature $\geq \bar{\bar \sigma}$ and geometry bounded by $\bar{\bar \cA}$, into the GL-sum 
			\begin{eqnarray*} 
				(B_\rho(p_1)  \,\#  &\!\!\!\! \underbrace{S^3   \, \#  \dots \, \# S^3 }\!\!\!\! &,\,  g \# \,  \can^{(n)} )\\
				&  n \, spheres & 
			\end{eqnarray*}
			where $\can^{(n)}$ denotes a metric obtained by GL-sum of $n$ round $3$-spheres of scalar curvature $\sigma$, made with parameters $(\rho,\sigma,\eta)$. Let $\widehat g_2^{\widehat V} \in \cR_{\sigma'}^{\cA'}(\widehat M)$ be the resulting metric. }\\
			
		\begin{proof}  Apply Lemma \ref{lem:conf-cyl} ($n-1$ times). \end{proof}
		
		 Next we choose points $x'(p_2), \ldots, x'(p_n) \in B_\rho(p_1) \# S^3(p_1)$ so that the balls $B_{6\rho}(x'(p_i))$ (for the metric  $g \# \,  \can^{(n)}$) are round, pairwise disjoint and at distance $\geq 6\rho'$ from $H(p_1)$ (this is possible since the number of $3$-spheres is large enough). Our next goal  is to deform $\widehat g_2^{\widehat V}$ by moving each base point $p_i$ to $x'(p_i)$, using paths in $V \# S^3(p_1)$. Before doing that, we need to adjust the parameters of the GL-sum to be compatible with the metric $g \# \can^{(n)}$ (a priori, the parameters $(\rho,\sigma,\eta )$ we begin with are compatible with $g$ only). Note that the metric $g \# \can^{(n)}$ has scalar curvature greater than $\sigma'$ and geometry bounded by $\cA'$. Hence it admits $(\rho',\sigma',\eta')$ as GL-parameters at any point, by definition of $\rho'$ and $\eta'$.

		Step 2. {\it For each $i \in \{2,\ldots,n\}$, deform $\widehat g_2^{\widehat V}$ on $B_\rho(p_i) \# S^3(p_i)$, without changing the metric near the boundary nor on the hemisphere $H(p_i)$, through metrics with scalar curvature greater than $\bar \sigma$ and geometry bounded by $\bar \cA$, into the GL-sum made at $\{(p_i,x(p_i))\}$ with the parameter $(\rho',\sigma',\eta')$. Let $\widehat g_3^{\widehat V} \in \cR_{\sigma''}^{\cA''}(\widehat M)$ be the resulting metric.}\\
		
		\begin{proof} Apply Lemma \ref{lem:conf-gamma}. \\ 
		\end{proof}

	    The crucial step is the following. \\
		
		Step 3. {\it There is an isotopy $(\widehat g_t^{\widehat V})_{t \in [3,4]}$ with scalar curvature greater than $\sigma''$ and geometry bounded by $\cA''$, which is constant on the hemispheres $H(p_i)$, $ i\geq 1$, such that $(\widehat V,\widehat g_4^{\widehat V})$ is isometric to the GL-sum 
			$$ \left(V\,\, \# \,\, S^3(p_1) \,\, \#\, \bigsqcup_{i \geq 2} S^3(p_i)\, ,\, g \,\, \# \,\, \can^{(n)}\,  \# \, \sqcup_{i\geq 2} \can\right) $$
			made at  $(p_1,x(p_1))$ and $\{(x'(p_i),x(p_i)) \mid i \geq 2 \}$  (with parameter $(\rho,\sigma,\eta)$ for $i=1$ and $(\rho',\sigma',\eta')$ for $i\geq 2$). }\\
		\begin{proof}
		Recall that $\can^{(n)}$ is a metric on $S^3$ isometric to a GL-sum of $n$ round spheres of scalar curvature $\sigma$. We will apply  Proposition \ref{prop:continuousGL}  iteratively on the manifold $M\, \# \, S^3(p_1) \sqcup \bigcup_{i \geq 2} S^3(p_i)$ with parameters $(\rho',\sigma',\eta')$, 
		moving points $p_i \leadsto x'(p_i)$ along paths in $V\, \# \, S^3(p_1)$. Proposition cannot be applied a priori by moving simultaneously the $(n-1)$ paths, because this would require having pairwise disjoint $3\rho'$-neighbourhoods for these paths. As the paths must traverse the neck  $B_\rho(p_1) \# B_\rho(x(p_1))$, this cannot be achieved with a radius $\rho'$ independent on $n$. The trick is to apply the proposition iteratively by moving only one point at a time, i.e. we move $p_i \leadsto x'(p_i)$, others points $x'(p_2),\ldots,x'(p_{i-1}), p_{i+1}, \ldots,p_n$ being fixed. In order to do this, it is sufficient to check the following claim:
		
		\paragraph{Claim}
		{\it For each $i \geq 2$, there exists a continuous path from $p_i$ to $x'(p_i)$ in 
			$$\left( V\, \# \, S^3(p_1) \right)\setminus \bigcup_{j \geq 2, j\not=i} B_{6\rho'}(p_j) \cup B_{6\rho'}(x'(p_j))$$
		and whose $3\rho'$-neighbourhood is contained in $V\, \# \, S^3(p_1) \setminus H(p_1)$}. 
		
Let prove this claim. Let $i \geq 2$. To begin we prove that there exists  a path in $U$ from $p_i$ to $p_1$, which remains at distance $\geq 6\rho'$ of all $p_j$ ($j\geq2, j \not=i$). By definition of $U$ (cf.~Lemma~\ref{lem:divide P}),  each ball $B_{6\rho'}(p_j)$ is compactly contained in $U$.  Note that $d(p_j,p_k) \geq 12\rho'$ for $k\not=j$, and that $6\rho'$ is smaller than the injectivity radius. Since $U$ is connected, there is  a path in $U$ from $p_i$ to $p_1$. If it intersects $B_{6\rho'}(p_j)$, we can replace the part of the path 
in  $B_{6\rho'}(p_j)$ by an arc in the boundary sphere  $S_{6\rho'}(p_j)$. This arc is contained in $U$ and at distance $\geq 6\rho'$ of any other $p_k$. Hence $U \setminus ( \cup_{j\geq2, j \not=i} B_{6\rho'}(p_j))$ is path-connected. 
Similarly in $(S^3(p_1), \can^{(n)})$, there exists a path joining $x(p_1)$ to $x'(p_i)$, disjoint from $B_{6\rho'}(x'(p_j))$ ($j\geq2, j \not=i$) and at distance $\geq 6\rho'$ from $H(p_1)$.  Using these two paths, we can find a path in the connected sum $U\, \# \, S^3(p_1)$ (made at $(p_1,x(p_1)$) joining $p_i$ to $x'(p_i)$ and at distance $\geq 6\rho'$ from $p_j$, $x'(p_j)$ and $H(p_1)$. Recalling that the $3\rho'$-neighbourhood of $U$ is contained in $V$, it follows that the $3\rho'$-neighbourhood of the path is contained in $V\, \# \, S^3(p_1) \setminus H(p_1)$.
This proves the claim, and concludes the proof of Step~3.

		\end{proof}

		Step 4. {\it For each $i \in \{2,\ldots,n\}$, deform  $\widehat g_4^{\widehat V}$ on $B_\rho(x'(p_i)) \# S^3(p_i)$, without changing the metric near the boundary nor on the hemisphere $H(p_i)$, through metrics with scalar curvature greater than $\bar \sigma$ and geometry bounded by $\bar \cA$, into the GL-sum made at $\{(x'(p_i),x(p_i))\}$ with the parameter $(\rho,\sigma,\eta)$. Let $\widehat g_5^{\widehat V} \in \cR_{\sigma'}^{\bar \cA'}(\widehat M)$ be the resulting metric.}\\
		
		\begin{proof}  Apply Lemma \ref{lem:conf-gamma} again. \\
(Note: There is an abuse of notation in the statement:  $\widehat g_4^{\widehat V}$ is a metric defined on $\widehat M$, and in particular on its submanifold  $B_\rho(p_i) \# S^3(p_i)$,  but can be seen 
as a metric on $B_\rho(x'(p_i)) \# S^3(p_i)$ only modulo  the isometry of Step 3.  ) \\
		\end{proof}
		
		 We define $(\widehat X_{\widehat V},g_{\widehat X_{\widehat V}})$ as the GL-sum 
	\begin{equation} 
 \left(S^3(p_1) \, \, \# \, \, S^3(p_2) \, \, \# \, \, \ldots \, \, \# \, \, S^3(p_n), \, \can^{(n)} \, \, \# \, \, \can \, \, \# \, \, \ldots \, \, \# \, \, \can \right) \label{eq:g_X}
\end{equation}
		made at $\{(x'(p_i),x(p_i)) \mid i \geq 2 \}$ with parameter $(\rho,\sigma,\eta)$.  We   have that 
		$(\widehat V,\widehat g_5^{\widehat V})$ is isometric to the GL-sum 
		$$\left( V \sqcup \widehat X_{\widehat V},  g  \sqcup g_{\widehat X_{\widehat V}}\right)_\#(\{(p_1,x(p_1)\})$$
		where the point $x(p_1) \in S^3(p_1)$ is considered as a point of $X_{\widehat V}$.  
		
		This concludes the construction of $\widehat g_t^{\widehat V}$ and $(\widehat X_{\widehat V},g_{\widehat X_{\widehat V}})$ when $N=1$.
		
Let us now deal with the general case. Write $[1,5]$ as an union of successive intervals, $[1,5]=I_1 \cup I_2 \cup \ldots \cup I_N$. Recall that $A=A_k=U_{i_1}  \cup \ldots \cup U_{i_N}$, and that $V$ is the $3\rho'$-neighbourhood of $A$.  On $I_1$ we use the four-step procedure above (modulo the obvious reparametrisation) on the set $U_{i_1}$. Let $t_1 = \max{I_1}$. The metric   $(\widehat V, \widehat g_{t_1}^{\widehat V})$ is isometric to the GL-sum of $(V,g)$ with  $(X_1,g_{X_1})$ (a finite connected GL-sum of round $3$-spheres as in \eqref{eq:g_X}) and a finite number of round $3$-spheres $S^3(p_\alpha)$. The points used for this GL-sum are the 
elements of $(\cP_{sep}\cap U_{i_1})  \cup (\cP \cap ( U_{i_2} \cup \ldots U_{i_N}))$ and corresponding points in $X_1$ and the $S^3(p_\alpha)$. Starting from this new GL-sum, apply on $U_{i_2}$  during the interval $I_2$ the four-step procedure as above. Then pull back the deformation to $\widehat V$ by the isometry. Let  $(X_2,g_{X_2})$  be the compact connected GL-sum  of round $3$-spheres created at this step. Iterate. The GL-manifold  $(\widehat X_{\widehat V},g_{\widehat X_{\widehat V}})$ is then defined as the disjoint union  $ \bigcup_{i=1,\ldots, N} (X_i,g_{X_i})$.

	    Having define isotopies $(\widehat g_t^{\widehat V_k})_{t \in [1,5]}$ on $\widehat M$, with support in $\widehat V_k$, we thus  define the isotopy $(\widehat g_t)_{t \in [1,5]}$ on $\widehat M$ by $$ \left\{\begin{array}{cccl}
		\widehat g_t & = & \widehat g_t^{\widehat V_k} & \textrm{ on } \widehat V_k,\, \forall k\\ 
		\widehat g_t & = & \widehat g_1 & \textrm{ on } \widehat M \setminus \bigcup_k \widehat V_k.
		\end{array}  \right.$$
	We define $(\widehat X,g_{\widehat X})$ as the disjoint union of the $(\widehat X_{\widehat V_k},g_{\widehat X_{\widehat V_k}})$. For each $p \in \cP_{sep}$ the point $x(p) \in S^3(p)$ can be considered as a point in $\widehat X$, and we set $\widehat x(p)=x(p)$. That finishes the proof of Sublemma~\ref{sublemma:deux}.
	
	\end{proof} 
		
We conclude the proof of Lemma \ref{lem:deform gl}. To define $\{g_\#^t\}_{t \in [1,5]}$, we reconnect $(\widehat M, \widehat g_t)$ into $(M_\#, g_\#^t)$ by doing the GL-sum 
$H(p^-) \# H(p^+)$ for all pairs of hemispheres, with parameters $(\rho,\sigma,\eta)$ (recall that $ \widehat g_t$ is constant on the hemispheres). Thus we set $g_\#^t=(g_t)_\#(\{(-x(p^-_\alpha),-x(p^+_\alpha))\}_\alpha)$. Similarly, we define $(X,g_X)$ as $(X,g_X)= (\widehat X,\widehat g_X)_\#(\{(-x(p^-_\alpha),-x(p^+_\alpha))\}_\alpha)$.	
\end{proof}

We conclude the proof of Proposition \ref{prop:sum_isotopies} using Lemmas \ref{lem:divide P} and \ \ref{lem:deform gl}. These lemmas imply that $\cP_{sep} \sqcup x(\cP_{sep}) \subset M \sqcup X$ is $3\rho$-separated for $g_t \sqcup g_X$. Therefore Proposition \ref{prop:continuousGL} applies on the manifold $(M \sqcup X,g_t \sqcup g_X)$ with basepoints $\{(p,x(p)) \mid p \in \cP_{sep}\})$ (which are fixed), and with parameters $\rho,\sigma,\eta$.  This proposition gives a deformation $(M \sqcup X, g_t \sqcup g_X)_\#$, with scalar curvature greater than $\sigma'$ and geometry bounded by $\cA'$, towards the GL-sum $(M \sqcup X, g_1 \sqcup g_X)_\#$, where $g_1$ and $g_X$ are GL-metrics. Let us denote by $\cS^M$ (resp. $\cS^X$) a 
		straight spherical splitting of $(M,g_1)$ (resp. $(X,g_X)$). After slightly deforming the GL-sum $(g_1 \sqcup g_X)_\#$ by moving base points $(p,x(p))$ (for $p \in \cP_{sep}$), we can assume that $\cS^M \sqcup \cS^X$ embeds in $(M \sqcup X)_\#$ and is straight for $(g_1 \sqcup g_X)_\#$. Denoting by $\cS^\#$ a collection of straight spheres associated 
		to the GL-sum at $\{(p,x(p)) \mid p \in \cP_{sep}\}$, then $\cS^\# \sqcup \cS^M \sqcup \cS^X$ is a straight spherical splitting of $(M \sqcup X, g_1 \sqcup g_X)_\#$. By applying dilations, we can finally arrange that the isotopy starting from $g_\#$ has scalar curvature greater than $9\sigma/10$, and geometry bounded by some $\bar \cA^{\#}$.

\end{proof}

 \section{Isotopy to a GL-metric: proof of Theorem \ref{thm:isotopy_to_GL}}\label{sec:isotopy_to_GL-can}
 
 The goal of this section is to prove  Theorem \ref{thm:isotopy_to_GL}, which we restate below.

\begin{theo}\label{thm:isotopy_to_GL-can} For every $A>0$ there exists $B=B(A)>0$ such that for every $3$-manifold $M$ and every metric $g \in \mathcal{R}_1^A(M)$, there exists a $GL$-metric $g' \in \mathcal{R}_1^{B}(M)$ isotopic to $g$ in 
$\mathcal{R}_1^{B}(M)$.
\end{theo}

\begin{proof}
Let $A>0$. Theorem~\ref{thm:existence_surgical} gives a number $\tau>0$ and a tuple $\cQ=\cQ_{[\epsi^{-1}]}$.

Let $M$ be a $3$-manifold and $g \in \cR_1^A(M)$. Fix an orientation of $M$. 
By Theorem~\ref{thm:existence_surgical} and Remark~\ref{rem:extinct}
we get  a surgical solution $(M(\cdot),g(\cdot))$ with geometry bounded by $Q_0$ and singular times $0<t_1 < \ldots < t_{j+1}<2$ such that $M(t)=\emptyset$ when $t \in (t_{j+1},2]$.  Thanks to Theorem \ref{thm:existence_surgical}, the number of surgeries $j+1$ is bounded by a number depending on $A$ only. Let us set $t_0=0$.

For convenience we set $M_0=M$ and $g_0(t)=g(t)$ on $[t_0,t_1]$. Moreover we set, for all $i\geq 1$, $M_i =M_+(t_i)$; we define $g_i$ on $[t_i,t_{i+1}]$ by $g_i(t_i)=g_+(t_i)$ and, for all  $t\in (t_i,t_{i+1}]$, $g_i(t)=g(t)$. 
We set $(M'_i, g'_i)=(M'(t_i), g'(t_i))$. Recall that $M'(t_i)$ is the union of the post-surgery manifold $M_i$ and all discarded components. Finally we set $M_{j+1}=\emptyset$.

Let $(\mathcal{P}_i)$ be the following property:

\medskip
 
 $(\mathcal{P}_i)$ : There exists $\cB^i=(B^i_k)_{0\leq k \leq [\epsi^{-1}]}$ and a map $h_i:[t_i, 2] \to \cR_1^{B_0^i}(M_i)$ such that
 \begin{enumerate}[\indent (1)] 
  \item $h_i(t_i)=g_i(t_i)$ and each component of $(M_i,h_i(2))$ is a GL-metric, 
  \item $h_i(t)$ has geometry bounded by $\cB^i$ when $t \in [\max(t_i,\tau),2]$. 
   \end{enumerate} 

We will prove that $(\mathcal{P}_i)$ holds for all $i\in\{0,\ldots,j\}$ by backward induction. This will show that $(\mathcal{P}_0)$ is true, which gives the required conclusion.
 
Let us first prove that $(\mathcal{P}_j)$ holds. We consider the Riemannian manifold $(M_j,g_j(t_j))$. We need to construct a metric isotopy $h_j(t)$ on $M_j$ for $t\in [t_j, 2]$.  

We set $h_j(t)=g_j(t)$ on $[t_j,t_{j+1}]$. Item 1) of Theorem \ref{thm:existence_surgical} shows that this path lies in $\cR_1^{Q_0}(M_j)$, and Item 2) of this theorem proves that $(2)$ of Property $(\mathcal{P}_j)$ holds in restriction to $[\max(t_j,\tau),t_{j+1}]$. We want to extend $h_j$  to the interval $[t_{j+1},2]$. Recall that the splitted  Riemannian manifold $(M'_{j+1},g'_{j+1})$ is $\epsi$-locally canonical since $M_{j+1}=\emptyset$. Lemma \ref{lem:loccan_to_GL} then applies to each component of $(M'_{j+1},g'_{j+1})$ and provides an isotopy $t \in [t_{j+1},2] \to g'_{j+1}(t) \in \cR_1^{\mathcal{B}'}(M'_{j+1})$, where $\cB'=\cB(\cQ)$, such that each component of $(M'_{j+1},g'_{j+1}(2))$ is a GL-metric. 
We denote by $(M'_{j+1})_\#$ the manifold obtained from $(M'_{j+1}, g'_{j+1})$ by GL-sum at tips of  pairs of caps produced by the metric surgery for appropriate parameters. We now identify $(M'_{j+1})_\#$ with $M_j$. 
Applying Proposition \ref{prop:sum_isotopies} and a rescaling, we deduce an isotopy $t \in [t_{j+1},2] \to h'_{j}(t) \in \cR_1^{\cB''}(M_j)$, starting from the GL-sum 
$((M'_{j+1})_\#,g'_{j+1}(t_{j+1})_\#)$, such that each component of $(M_j,h'_{j}(2))$ is a GL-metric. Using Lemma \ref{lem:gl-surg} in a small interval near $t_{j+1}$ in $[t_{j+1},2]$, we isotope each surgery neck of $(M_j,h_j(t_{j+1}))$ to the corresponding GL-sum of caps of $(M_j,h'_j(t_{j+1}))$. Denoting $h_j$ the corresponding isotopy on $[t_{j+1},2]$, we obtain a continuous path $h_j$ on $[t_j,2]$. This proves Property~$(\mathcal{P}_j)$.
 
Let us now prove that $(\mathcal{P}_{i+1}) \Rightarrow (\mathcal{P}_{i})$. We assume that $(\mathcal{P}_{i+1})$ holds.
Then there exists $\cB^{i+1}=(B^{i+1}_k)_{0\leq k \leq [\epsi^{-1}]}$ and a  continuous path $h_{i+1}:[t_{i+1},2] \to \cR_1^{B_0^{i+1}}(M_{i+1})$  satisfying (1) and (2). 
 We have to prove that there exists $\cB^i$ such that each connected component of $(M_{i},g_{i}(t_{i}))$ is isotopic in $\cR_1^{B_0^{i}}$ to a GL-metric by an isotopy satisfying (2). 
 On $[t_i,t_{i+1}]$ we set $h_i(t) = g_i(t)$ as before. Then 
$M'_{i+1} \supset M_{i+1}$ and $M'_{i+1} \setminus M_{i+1}$ is the union of the discarded components. The components of $(M'_{i+1},g'_{i+1})$ are now of two types:  those of $M_{i+1}$, which satisfy $(\mathcal{P}_{i+1})$ by the induction assumption, and 
 those of $M'_{i+1} \setminus M_{i+1}$, which are $\epsi$-locally canonical for $g'_{i+1}$. On $M_{i+1}$ we denote by $g'_{i+1}(t)$ an isotopy defined for $t\in [t_{i+1}, 2]$ and given by $(\mathcal{P}_{i+1})$.
 On $M'_{i+1} \setminus M_{i+1}$ we denote by $g'_{i+1}(t)$ the isotopy defined for $t\in [t_{i+1}, 2]$ and given by Lemma \ref{lem:loccan_to_GL}. As above we obtain an isotopy $h'_i(t)$ on $M_i$ defined for $t\in [t_{i+1}, 2]$, by applying  Proposition \ref{prop:sum_isotopies} and a rescaling to  $(M'_{j+1},g'_{j+1})$. We conclude again using Lemma \ref{lem:gl-surg}.

This completes the proof of Theorem~\ref{thm:isotopy_to_GL-can}.  
\end{proof}

\end{spacing}

\appendix

\section{Construction of metrics with uniformly positive scalar curvature}\label{sec:appendix}

We prove Theorem~\ref{thm:BBM11}, which we restate for the convenience of the reader.
\begin{theo}
Let $M$ be an oriented, connected $3$-manifold. Then $M$ admits a complete Riemannian metric of uniformly positive scalar curvature and bounded geometry if and only if there exists a finite collection $\cF$ of spherical $3$-manifolds such that $M$ is a connected sum of members of $\cF$.
\end{theo}

\begin{proof}
The `only if' part is proven in~\cite{B2M:scalar}, except that the conclusion there allows for factors which are diffeomorphic to $S^2\times S^1$. Those factors can be removed by adding extra edges in the graph.

We prove the `if' part using the material of Section~\ref{sec:GL}.
Let $(G,X)$ be a pair presenting $M$ as a connected sum. First we notice that the graph can be modified so that every vertex has degree at most three: for each vertex $v$ of degree $d\ge 4$ (if any), replace $v$ by a finite tree $T_v$ with $d$ leaves and such that every vertex of $T_v$ has degree at most 3; for each such $v$, fix an arbitrary vertex $w_v$ of $T_v$; then associate $X_v$ to $w_v$ and $S^3$ to the other vertices of $T_v$. The resulting manifold is diffeomorphic to $M$.

Next we put on each $X_v$ an arbitrary Riemannian metric with scalar curvature $\ge 1$, geometry bounded by some constant $C$ and such that each $X_v$ contains a disjoint union of three metric balls of radius $C^{-1/2}$. We fix a number $\rho\in (0,C^{-1/2})$. In each $X_v$ we choose a finite set of points $x_i$ with cardinality the degree of $v$, and such that the metric balls of radius $\rho$ around the $x_i$'s are pairwise disjoint. We remove those metric balls, getting a collection of punctured $3$-manifolds. Let $Y$ be their disjoint union. Finally, we choose a neck $N_e$ for each edge of $G$ and glue $N_e$ to the corresponding boundary spheres of $Y$. The resulting manifold is diffeomorphic to $M$; by our generalisation of the Gromov-Lawson construction (see Subsection~\ref{subsec:gl}), it carries a complete metric with uniformly positive scalar curvature and bounded geometry.
\end{proof}

\section{Deforming metrics of positive scalar curvature on closed manifolds}\label{appendix:marques}
We explain how to prove Theorem~\ref{thm:marques deforming2}, which we restate for convenience:
\begin{theo}\label{thm:marques deforming3}
Let $M$ be a closed, oriented 3-manifold such that $\mathcal{R}_1(M)\not=\emptyset$. Then $\mathcal{R}_1(M)/\diff^+(M)$ is path-connected  in the $\mathcal{C}^\infty$ topology.  
\end{theo}

\begin{proof}
There are three differences between Theorem~\ref{thm:marques deforming} and Theorem~\ref{thm:marques deforming2}: (1) the manifold $M$ need not be connected; (2) we work with $\cR_1(M)$ as opposed to $\cR_+(M)$; (3) we work with $\diff^+(M)$ instead of $\diff(M)$.

Point (1) does not pose any problem since we can work with each component separately. Point (2) is dealt with by the following trick:

\begin{lem}
Let $g,g'$ be two metrics in $\cR_1(M)$ which are isotopic in $\cR_+(M)$. Then $g,g'$ are isotopic in $\cR_1(M)$.
\end{lem}

\begin{proof}
Let $g_t$ be an isotopy from $g$ to $g'$. Since $M$ is compact, there exists $\sigma>0$ such that $g_t\in\cR_\sigma$ for all $t$. Thus Lemma~\ref{lem:sigma} applies.
\end{proof}

To deal with Point~(3), we need to track down all places in the article~\cite{marques:deforming} where diffeomorphisms are introduced. Recall that the proof has two parts: show that any $g \in\cR_+(M)$ is isotopic to a canonical metric; show that canonical metrics on $M$ are isotopic modulo diffeomorphism. Canonical metrics are defined on pages~841 and~842: a metric $g$ on $M$ is canonical if it is isometric to $(\widehat M,\widehat g)$, a GL-sum of a standard $S^3 \subset \RR^4$ and of spherical manifolds with constant sectional curvature metrics (for some choice of base points and orthonormal bases). We can always assume this isometry to be positive, orienting $\widehat M$ by the choice of the bases, reversing all bases if necessary. The important point is that arguments on page~842, based on Milnor's and De Rham's theorems, yield that any two canonical metrics on $M$ lie in the same path-connected component of the moduli space $\cR_+(M)/\diff^+(M)$. 
For the second part of the proof, it is thus unnecessary to control the sign of diffeomorphisms, if any. 
\end{proof}

\bibliographystyle{alpha}
\bibliography{biblio}

\Addresses

\end{document}